\documentclass[a4paper,11pt]{article}
\normalsize

\usepackage{etoolbox}
\makeatletter
\renewcommand*\l@section{\@dottedtocline{1}{1.5em}{2.3em}}
\patchcmd{\@dottedtocline}{\leavevmode}{\leavevmode\bfseries\boldmath}{}{}
\patchcmd{\@dottedtocline}{\normalfont}{\normalfont\bfseries\boldmath}{}{}
\patchcmd{\l@part}{\bfseries}{\fbfseries\boldmath}{}{}
\patchcmd{\l@section}{\bfseries}{\bfseries\boldmath}{}{}
\patchcmd{\@part}{\bfseries}{\bfseries\boldmath}{}{}
\patchcmd{\@spart}{\bfseries}{\bfseries\boldmath}{}{}
\patchcmd{\section}{\bfseries}{\bfseries\boldmath}{}{}
\patchcmd{\subsection}{\bfseries}{\bfseries\boldmath}{}{}
\patchcmd{\subsubsection}{\bfseries}{\bfseries\boldmath}{}{}
\patchcmd{\paragraph}{\bfseries}{\bfseries\boldmath}{}{}
\patchcmd{\subparagraph}{\bfseries}{\bfseries\boldmath}{}{}
\makeatother

\usepackage[dvipsnames]{xcolor}
\definecolor{myred}{RGB}{211, 32, 32}
\definecolor{myblue}{RGB}{0, 76, 153}
\definecolor{mygreen}{RGB}{0, 180, 0}
\definecolor{mypurp}{RGB}{153, 0, 153}

\usepackage[style=alphabetic,maxbibnames=99]{biblatex}
\usepackage[affil-it]{authblk}

\usepackage{mathtools}
\usepackage{amsfonts, amsmath, amssymb}
\usepackage{bm}
\usepackage{upgreek}
\usepackage[mathscr]{euscript}
\usepackage{physics}
\usepackage{dsfont}
\usepackage{tabularx}
\usepackage{faktor}
\usepackage{dsfont}
\usepackage{amsthm}
\usepackage{resizegather}

\makeatletter
\newcommand*{\rom}[1]{\expandafter\@slowromancap\romannumeral #1@}
\makeatother

\let\oldsection\section
\renewcommand{\section}{
  \renewcommand{\theequation}{\thesection.\arabic{equation}}
  \oldsection}
\let\oldsubsection\subsection
\renewcommand{\subsection}{
  \renewcommand{\theequation}{\thesubsection.\arabic{equation}}
  \oldsubsection}

\DeclareMathOperator{\End}{End}
\DeclareMathOperator{\Hom}{Hom}

\newcommand{\I}{\tilde{I}}
\newcommand{\kk}{\widetilde{\mathsf{k}}}

\newcommand{\ldual}[1]{\prescript{L}{}{{{#1}}}}
\newcommand{\hf}[1]{\widehat{\mathfrak{#1}}}
\newcommand{\ch}[1]{\check{#1}}

\newcommand{\defeq}{\vcentcolon=}
\newcommand{\diff}{\text{d}}
\newcommand{\res}[1]{\text{res}_{#1}}

\newcommand{\ZZ}{\mathbb{Z}}
\newcommand{\CC}{\mathbb{C}}
\newcommand{\PP}{\mathbb{P}}
\newcommand{\g}{\mathfrak g}
\newcommand{\V}{\mathcal V}
\newcommand{\VV}{\mathbb V}
\newcommand{\der}[1]{[#1]}
\renewcommand{\sl}{\mathfrak{sl}}
\newcommand{\vak}{\ket{0}^{k}}
\newcommand{\vakk}{\ket{0}^{\bm k}}
\newcommand{\vac}{\ket{0}}

\newtheorem{theorem}{Theorem}[subsection]
\newtheorem{proposition}{Proposition}[subsection]
\newtheorem{corollary}{Corollary}[subsection]
\theoremstyle{definition}

\newtheorem*{remark*}{Remark}
\allowdisplaybreaks

\newenvironment{mproof}[1][\proofname]{%
\vspace{.4cm}\hrule
  \begin{proof}
}{%
  \end{proof}
  \hrule\vspace{.4cm}
}

\usepackage{enumitem}
\usepackage[
    colorlinks = true,
    linkcolor=NavyBlue,          
    citecolor=BrickRed,        
    filecolor=magenta,      
  unicode,
  ]{hyperref}
 \def\eqref#1{(\textcolor{NavyBlue}{\ref{#1}})}
  \usepackage{url}
\usepackage{cleveref}

\usepackage{geometry}
\usepackage{authblk}
\usepackage{titletoc}
\allowdisplaybreaks

\usepackage{appendix}

\title{\sffamily \bfseries Quartic Hamiltonians, and higher Hamiltonians at next-to-leading order, for the affine $\bm{\sl_2}$ Gaudin model}
\author{Tommaso Franzini\footnote{email \texttt{\href{mailto:t.franzini@herts.ac.uk}{t.franzini@herts.ac.uk}}} $\qquad$ Charles Young\footnote{email \texttt{\href{mailto:c.young8@herts.ac.uk}{c.young8@herts.ac.uk}}}}
\affil{\small Department of Physics, Astronomy and Mathematics,\\ University of Hertfordshire, Hatfield, AL10 9AB, United Kingdom}

\date{}

\begin{document}
\maketitle
\thispagestyle{empty}
\begin{abstract}
In this work we will use a general procedure to construct higher local Hamiltonians for the affine $\sl_2$ Gaudin model.
We focus on the first non-trivial example, the quartic Hamiltonians. We show by direct calculation that the quartic Hamiltonians commute amongst themselves and with the quadratic Hamiltonians which define the model.

We go on to introduce a certain next-to-leading-order semi-classical limit of the model. In this limit, we are able to write down the full hierarchy of higher local Hamiltonians and prove that they commute.
\end{abstract}

\setcounter{tocdepth}{1}
\tableofcontents
\noindent\rule{\textwidth}{1pt}
\newpage

\section{Introduction}
Gaudin models represent interesting theories that find applications in several contexts in mathematical physics. In particular, they provide a general framework to study rich classes of classical and quantum integrable 1+1 dimensional integrable theories \cite{Gaudin1,Gaudin2,Vicedo2018,Lacroixthesis}. Recently, it has also been found that classical untwisted Gaudin models provide a dual description of such theories to the one given by 4D holomorphic Chern-Simons theories \cite{V2021,LV2021}.

In this paper we will focus on \emph{quantum Gaudin models}.
To define them, one needs a collection of distinct points $\{z_1,\dots,z_N\}\in\CC\PP^1$ on the Riemann sphere and a Kac-Moody algebra $\g$, which can be of finite or affine type.
These models are defined by the quadratic Hamiltonians, namely
\begin{equation}\label{hamiltonian}
    \mathscr{H}_i=\sum_{\substack{j=1 \\ j\neq i}}^N\kappa_{ab} \frac{I^{a (i)}I^{b(j)}}{z_i-z_j}, \qquad i=1,\dots,N,
\end{equation}
in the (suitably completed, in the affine case) tensor product $U(\g)^{\otimes N}$. Here $I^{a (i)}$ are the generators of $\g$, defined at site $i$, and $\kappa$ is the invariant bilinear form on $\g$.

Finite-type quantum Gaudin models have been deeply studied. In particular, it is known that the quadratic Hamiltonians \eqref{hamiltonian} are part of a larger family of mutually commuting operators $\mathscr{B}\in U(\g^{\otimes N})$ called the Bethe or Gaudin subalgebra  \cite{Frenkel2004,Talalaev,RYB2006}. Moreover, these higher Hamiltonians can be diagonalized with an elegant form of Bethe Ansatz, where the eigenvector is the Schechtman-Varchenko vector and the eigenvalues are encoded as functions on a space of $\ldual\g$-opers, the Langlands dual algebra of $\g$  \cite{FFR,MTV2006,MTV2007,MV2002,MV2004}.

On the other hand, affine-type Gaudin models are still far from being fully understood.

One approach to affine Gaudin models was proposed in a recent paper \cite{KL2021}. In this work the authors propose an integrable model called \emph{generalized} affine $\mathfrak{sl}_2$ Gaudin model, which reproduces the usual $\hf{sl}_2$ model in a certain limit. Their construction is based on a new realization of the subalgebra $U_q(\mathfrak b_-)$ through a new class of vertex operators, and fits affine Gaudin models into the general procedure first given in the seminal work \cite{BLZ1}.

\begin{figure}[h!]
    \centering
\includegraphics[scale=.7]{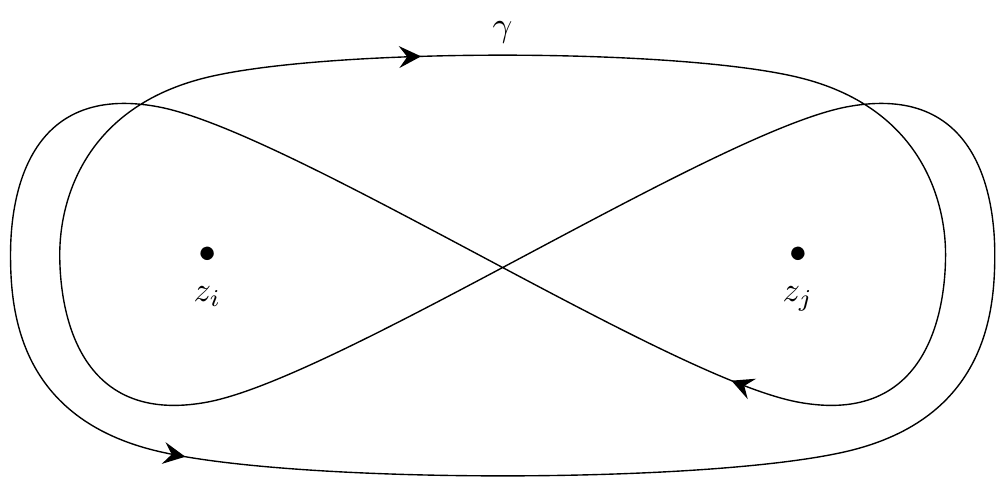}
    \caption{An example of Pochhammer contour $\gamma$ around any two marked points.}
    \label{fig:pochhammer}
\end{figure}

Another approach, first proposed in the pioneering work \cite{FF2007}, is to treat affine Gaudin models in closer analogy with their finite-type cousins. In particular, the spectrum of higher Hamiltonians should be described by suitable functions on a space of opers. Some further conjectures of how this might work, at least for the local Hamiltonians, were made in \cite{LVYaffine}, where it was conjectured that the eigenvalues of higher local Hamiltonians of the affine Gaudin models, as well as the Hamiltonians themselves, are given by hypergeometric-type integrals in the spectral plane, namely
\begin{equation}
    \widehat{Q}_n^{\gamma} = \int_{\gamma} \mathscr{P}(z)^{-n/2} \varsigma_n(z)_{[0]} \diff z ,
\end{equation}
where $n$ lives is a (multi)set of indices given by the exponents of the algebra $\g$, $\mathscr{P}$ is a certain multi-valued function defined by the data of the levels $k_i$ of the modules attached to the marked points $z_i$, $\gamma$ is a Pochhammer contour in the spectral plane around any two of these points (see e.g. \cref{fig:pochhammer}) and $\varsigma_n(z)_{[0]}$ can be thought as the Hamiltonian density.

The key step in computing the higher Hamiltonians is to characterize these Hamiltonian densities, which are obtained by defining a suitable state $\varsigma_n(z)\in\mathcal V$ for each given exponent $n$. In order to do that, it is possible to exploit the general properties consistent Hamiltonians must obey: they have to commute with the generators $\{I^{\alpha}_n\}_{n\in\ZZ}^{\alpha=1,\dots,\dim\mathring{\g}}$ of the algebra $\g$ defining the model as well as amongst themselves (here $\mathring{\g}$ denotes the underlying finite algebra). As we will see in \cref{commutingham}, this is equivalent to the following requirements
\begin{equation}\begin{split}\label{tsing}
    \Delta I^{\alpha}_{n\geq 0} \varsigma_m(z) = 0 &\mod \text{twisted derivatives}, \\
    \varsigma_n(z)_{(0)}\varsigma_m(z) = 0 &\mod \text{twisted derivatives and translates},
\end{split}
\end{equation}
where $\Delta I^{\alpha}_{n\geq 0}$ represents the diagonal action of the positive modes of the generators of the algebra $\g$ and the zero mode $\varsigma_n(z)_{(0)}$ is intended in the vertex algebra sense (see \cref{VAstructure}).
We give a precise definition of translation and twisted derivative of a state in \cref{VAstructure,twistedder}.

The general expectation is that there exists a state $\varsigma_m(z)$ for every exponent $m$ of the affine algebra $\g$, and that it takes the following form
\begin{equation}\label{introtop}
    \varsigma_m(z) = t_{i_{1},\dots,i_{m+1}}I^{i_1}_{-1}(z)\cdots I^{i_{m+1}}_{-1}(z)\vac + \text{quantum corrections},
\end{equation}
where $I(z)=\sum_iI^{(i)}/(z-z_i)$ and $t$ is a certain symmetric invariant tensor of $\mathring{\g}$.
This particular structure is justified by the semi-classical counterpart of these models, which have been thoroughly studied \cite{Evans_1999,Evans_2001,LMV2017}; in particular, the precise choice of symmetric invariant tensors needed to ensure that the Hamiltonians Poisson-commute is well understood and related to Drinfel'd-Sokolov reduction \cite{Evans_2001,Drinfeld1985}. In the very simplest cases, including the cubic Hamiltonian in type $\widehat{\mathfrak{sl}}_{M\geq3}$, there are no quantum corrections needed \cite{LVYcubic}. The first case in which quantum corrections are present appeared in the recent paper \cite{LKT22}: it is the example of the quartic Hamiltonian in type $\widehat{\mathfrak{sl}}_2$. In this paper, we consider this example in detail. In particular, we carry out the full computation to show that the resulting quartic Hamiltonians commute amongst themselves (see \cref{th2} and \cref{gaudinham}).

This direct calculation reveals the following fact: it turns out that \emph{any} density $\varsigma_3(z)$ obeying the first of the two conditions in \cref{tsing} (i.e. the one which ensures the Hamiltonian commutes with the generators of the affine algebra), automatically also obeys the second condition (which is needed for the Hamiltonians to commute amongst themselves) at least for $n=1,3$, as shown in \cref{cor} below. We should stress that this property is, for the moment, highly non-obvious and suggestive; it would be good to get a more systematic understanding of why it should be true.

\bigskip

Already in this case of the exponent $n=3$, i.e. of quartic Hamiltonians, the direct computations needed are very lengthy. This is especially true of the computations needed to show the mutual commutativity of the quartic Hamilonians. For higher exponents $n \geq 5$, direct calculations become computationally difficult even with the aid of computer algebra, but we are able to prove a result for all exponents by truncating to the next-to-leading order in $\hbar$.
To introduce the dependence from the formal parameter $\hbar$, we perform a re-scaling of the generators, namely $I\to\tilde{I} \equiv \hbar I$, $k\to\tilde{k} \equiv k$, in such a way that every time we perform a commutation, we introduce a factor of $\hbar$. This procedure allows us to identify different quantum corrections by their $\hbar$ dependence.
In this spirit, we will then work modulo terms of order $\hbar^2$.
We show that, modulo such terms, the Hamiltonian density for each exponent $2n-1$, $n\in \mathbb{Z}_{\geq 1}$, of $\widehat{\mathfrak{sl}}_2$, takes the following form
\begin{equation*}\begin{split}
    \tilde{\varsigma}&_{2n-1}(z)= t_{i_1,\dots,i_{2n}} \I^{i_1}_{-1}(z)\I^{i_2}_{-1}(z)\cdots \I^{i_{2n}}_{-1}(z)\vac\\
    &\qquad+ \hbar \frac{n(2n+1)(2n-2)}{(2n-1)} t_{i_1,\dots,i_{2n-4}}f^{abc} \I^{a}_{-2}(z)\I^{b\prime}_{-1}(z)\I^{(c}_{-1}(z)\I^{i_1}_{-1}(z)\cdots \I^{i_{2n-4})}_{-1}(z)\vac
\end{split}
\end{equation*}
and prove that the resulting Hamiltonians commute up to and including terms of order $\hbar^2$.
\begin{center}
    ***
  \end{center}

\newpage
The paper is organized as follows.

In \cref{sec:VVmod} we recall the main ideas behind the theory of Kac-Moody algebras, their local completion and the concept of vacuum Verma module. We also recall the definition and the main features of vertex algebras.

In \cref{Construction} we define the algebra of observables of the Gaudin model and we recall the definition of invariant tensor of an algebra $\g$, focusing in particular on the case of $\hf{sl}_2$.

In \cref{quartic} we define the meromorphic states $I(z)$, which represent the building blocks for the states we want to define. We also present their commutation relations and we describe the gradation they give rise, which will be of fundamental importance to prove the main results of the paper.
After that, we characterize the space of states that vanish under the diagonal action of zero modes: this is the first step to define higher Gaudin Hamiltonians. At this point we will present all quantum corrections to the quartic state and we will explain why they appear in this specific example and not in the other known ones.
In the following subsection, we restrict the number of possible states by asking they be singular up to twisted derivative under the diagonal action of positive modes. We focus on the quadratic and quartic states.
Afterwards, we prove the main result of the paper, where we also ask that the 0th product (in the vertex algebra sense) of these states vanishes modulo twisted derivatives and translates. We will show that there is one quartic state, up to re-scaling and modulo the addition of translates and twisted derivatives, that satisfies this last requirement.
At this point, having a precise definition of the quartic state, we describe the general construction of the higher Hamiltonians, in the spirit of \cite{LVYaffine}. We will explain in detail why we ask for these properties and why they are important at the level of the quantum Hamiltonians.

In \cref{subleading} we try to solve the same problem from a different point of view.
Instead of focusing on one specific Hamiltonian, trying to work out its explicit definition, we want to characterize all Hamiltonians at arbitrary $n$, but working at sub-leading order. In order to do that we introduce the formal parameter $\hbar$ by making a particular re-scaling of the generators of the algebra. We will then prove similar theorems to those of the previous sections, working modulo terms at order $\hbar^3$.

\bigskip

\noindent\textbf{Acknowledgements}. The authors would like to thank Sylvain Lacroix for helpful discussions and for providing many detailed comments.
The research of Charles Young is supported by the Leverhulme Trust, Research Project Grant number RPG-2021-092.

\section{Vacuum Verma modules for \texorpdfstring{$\hf{sl}_2$}{affine sl2}}\label{sec:VVmod}
\subsection{Loop realization of \texorpdfstring{$\hf{sl}_2$}{affine sl2}}\label{subsec:loopreal}

We define the loop algebra $\sl_2[t,t^{-1}] = \sl_2\otimes\CC[t,t^{-1}]$ as the algebra of Laurent polynomials in a formal variable $t$ with coefficient in the finite-dimensional Lie algebra $\sl_2$. The Lie brackets on this algebra are given by
\begin{equation}\label{affinecommrel}
    [a\otimes f(t),b\otimes g(t)] = [a,b]_{\sl_2} \otimes f(t)g(t),
\end{equation}
where $f(t)$ and $g(t)$ are arbitrary Laurent polynomials in $\CC[t,t^{-1}]$.

Let $(\cdot\vert\cdot):\sl_2\times\sl_2\rightarrow\CC$ be the canonically normalized symmetric invariant bilinear form on $\sl_2$. It is given by taking the trace in the defining representation:
\begin{equation}\label{killing}
    (a\vert b) \defeq \Tr(ab),
\end{equation}

It is possible to extend the loop algebra by a one-dimensional central element $\CC\mathsf{k}$,
\begin{equation}
    0 \longrightarrow \CC\mathsf{k} \longrightarrow \hf{sl}_2 \longrightarrow \sl_2[t,t^{-1}] \longrightarrow 0.
\end{equation}
This extension is called the affine Lie algebra $\hf{sl}_2$, whose commutation relations are
\begin{gather}
      [a\otimes f(t),b\otimes g(t)] = [a,b]_{\sl_2} \otimes f(t)g(t) - ( \res{t} f \diff g) (a\vert b)\mathsf{k},\\
    [\mathsf{k},\cdot] = 0.
\end{gather}
We shall use the notation
\begin{equation}
a_n \defeq a\otimes t^{n},\quad \text{for $a\in\sl_2$ and $n\in\ZZ$}.\label{andef}
\end{equation} The commutation relations can be equivalently written as
\begin{equation}\label{affinecommrel2}
    [a_m,b_n] = [a,b]_{n+m} - n \delta_{n+m,0}(a\vert b) \mathsf{k}.
\end{equation}
We can add to this algebra a one-dimensional derivation $\mathsf{d}$, such that $[\mathsf{d},\mathsf{k}] = 0$ and $[\mathsf{d},a\otimes f(t)] = a\otimes t\partial_tf(t)$, for all $a\in\sl_2$ and $f(t)\in\CC[t,t^{-1}]$. It is possible to show that this algebra is isomorphic to the Kac-Moody algebra over $\CC$ of type $\mathsf{A}_{1}^{(1)}$, see e.g. \cite[ch. 7]{kacinfinite},
\begin{equation}
    \mathfrak{g} = \hf{sl}_2\oplus\CC\mathsf{d}.
\end{equation}

\subsection{\texorpdfstring{$\hf{sl}_2$}{TEXT} as a Kac-Moody algebra}\label{KMdata}
The Cartan matrix for the Kac-Moody algebra of type $^{1}\mathsf{A}_{1}$ is defined as $A = (a_{i,j})_{i,j=0}^1 = (2\delta_{i,j}-\delta_{i+1,j}-\delta_{i-1,j})_{i,j=0}^1$. The Cartan decomposition is given by $\mathfrak{g}=\mathfrak{n}_{-}\oplus\mathfrak{h}\oplus\mathfrak{n}_{+}$. The Chevalley-Serre generators are $\{e_i\}_{i=0}^1\subset \mathfrak{n}_{+}$, $\{f_i\}_{i=0}^1\subset\mathfrak{n}_{-}$ while $\{\ch{\alpha}_i\}_{i=0}^1\subset\mathfrak{h}$ and $\{\alpha_i\}_{i=0}^1\subset\mathfrak{h^{\ast}}$ are respectively a basis for the Cartan subalgebra of simple coroots of $\mathfrak{g}$ and a basis for the dual Cartan subalgebra of simple roots of $\mathfrak{g}$. The latter are related via the canonical pairing between the Cartan algebra and its dual,  $\langle\cdot,\cdot\rangle:\mathfrak{h}^{\ast}\times\mathfrak{h}\rightarrow \CC$
\begin{equation}
    \langle\alpha_i,\ch{\alpha}_j\rangle=a_{i,j}.
\end{equation}
The fundamental commutation relations in $\mathfrak{g}$ are
\begin{equation}
\begin{gathered}
    [x,e_i]= \langle\alpha_i,x\rangle e_i,     \qquad      [x,f_i]= -\langle\alpha_i,x\rangle f_i,\\
     [x,x']= 0,        \qquad    [e_i,f_j]= \ch{\alpha}_i\delta_{ij},
\end{gathered}
\end{equation}
where $x,x'\in \mathfrak{h}$ and $i,j=0,1$, together with the Serre relations
\begin{equation}
    (\text{ad}e_i)^{1-a_{ij}}e_j=0, \qquad (\text{ad}f_i)^{1-a_{ij}}f_j=0.
\end{equation}

The Kac-Moody algebra $\mathfrak{g}$ has a central element $\mathsf{k}=\sum_{i=0}^1\ch{\alpha}_i$, which spans a one-dimensional centre.
A basis for the Cartan subalgebra $\mathfrak{h}$ is given by the coroots $\{\ch{\alpha}_i\}_{i=0}^1$ together with the derivation element $\mathsf{d}$, which by definition satisfies
\begin{equation}
    \langle\alpha_i,\mathsf{d}\rangle = \delta_{i,0}.
\end{equation}
If we remove the 0th row and column from $A$, we obtain the Cartan matrix corresponding to the finite dimensional Lie algebra $\sl_2$. This subalgebra of $\mathfrak{g}$ is generated by $e_1\in\mathfrak{n}_+$, $f_1\in\mathfrak{n}_-$ and $\ch{\alpha}_1\in\mathfrak{h}$.

\subsection{Local completion and vacuum Verma modules}\label{LCandVVM}
For any $k\in\CC$, let us define $U_k(\hf{sl}_2)$ as the quotient of the universal enveloping algebra $U(\hf{sl}_2)$ of $\hf{sl}_2$ by the two-sided ideal generated by $\mathsf{k}-k$. For each $n\in \ZZ_{\geq 0}$, let us introduce the left ideal $J_n = U_k(\hf{sl}_2)\cdot (\sl_2\otimes t^n\CC[t])$.
The inverse limit
\begin{equation}
    \widetilde{U_k}(\hf{sl}_2) = \varprojlim \faktor{U_k(\hf{sl}_2)}{J_n}
\end{equation}
is a complete topological algebra, called the local completion of $U_k(\hf{sl}_2)$ at level $k$.
With this definition, the elements of $\widetilde{U_k}(\hf{sl}_2)$ are possibly infinite sums of the type $\sum_{m\geq0} a_{m}$ of elements $a_m\in U_k(\hf{sl}_2)$ which do truncate to finite sums when working modulo any $J_n$.

A module $\mathscr{M}$ over $\hf{sl}_2$ is said to be smooth if, for all $a\in\sl_2$ and all $v\in\mathscr{M}$, $a_n v =0$ for sufficiently large $n$. A module $\mathscr{M}$ has level $k$ if $\mathsf{k}-k$ acts as zero on $\mathscr{M}$. Any smooth module of level $k$ over $\hf{sl}_2$ is also a module over the completion  $\widetilde{U_k}(\hf{sl}_2)$.

We can identify the subalgebra of positive modes $\sl_2[t]\oplus \CC\mathsf{k}\subset \hf{sl}_2$ and introduce the one-dimensional representation $\CC\vak$ defined by
\begin{equation}
    (\mathsf{k}-k)\vak=0, \qquad\qquad a_n\vak=0 \quad \text{for all $n\geq0$, $a\in\sl_2$}.
\end{equation}

We define $\mathbb{V}_0^k$, the vacuum Verma module at level $k$, as the induced smooth $\hf{sl}_2$-module
\begin{equation}
    \mathbb{V}_0^k = U(\hf{sl}_2)\otimes_{U(\sl_2[t]\oplus \CC\mathsf{k})}\CC\vak
\end{equation}
This vector space is spanned by monomials of the form $a_p\cdots b_q\vak$, with $a,\dots,b\in\sl_2$ and strictly negative mode numbers $p,\dots,q\in\ZZ_{<0}$. We call these vectors states.

Let us denote by $[T,\cdot]$ the derivation on $U_k(\hf{sl}_2)$ defined by $[T,a_n]=-na_{n-1}$ and $[T,1]=0$.
By setting $T(X\vak)=[T,X]\vak$ for any $X\in U_k(\hf{sl}_2)$, one can interpret $T$ as a translation operator $T: \mathbb{V}_0^k\rightarrow \mathbb{V}_0^k$: the reason for this identification will be clear in the next section.

\subsection{Vertex algebra structure}\label{VAstructure}
As we now recall, the vacuum Verma module defined in the previous section has the structure of vertex algebra.
Namely, we have the state-field map $Y(\cdot,x)$, which for every state $A\in\mathbb{V}_0^k$ associates a formal power series in the variable $x$,
    \begin{equation}\label{ymap}
    \begin{split}
        Y(\cdot,x): \mathbb{V}_0^k &\longrightarrow \Hom(\mathbb{V}_0^k, \mathbb{V}_0^k ((x)))\\
                    A &\longmapsto Y(A,x) = \sum_{n\in\ZZ}A_{(n)}x^{-n-1}
    \end{split}
    \end{equation}
where $A_{(n)}\in \End(\mathbb{V})$ is the $n$\textsuperscript{th} mode of $A$.
By definition, if $A=a_{-1}\vak$ for some $a\in\sl_2$, then $A_{(n)}=a_n$ for all $n\in\ZZ$, \textit{i.e.}
\begin{equation}
    Y(a_{-1}\vak,x)=\sum_{n\in\ZZ}a_n x^{-n-1}.
\end{equation}
The fields for all other states can be obtained with the following properties
\begin{equation}\label{recursiveY}
    Y(TA,x)=\partial_z Y(A,x), \qquad\qquad Y(A_{(-1)}B,x)=:Y(A,x)Y(B,x):
\end{equation}
where we have introduced the normal ordered product between fields
\begin{equation}\label{normalord}
        :Y(A,x)Y(B,x):\, = \left(\sum_{m<0}A_{(m)}z^{-m-1}\right)Y(B,x) +Y(B,x) \left(\sum_{m\geq0}A_{(m)}x^{-m-1}\right).
\end{equation}
In fact, any state $C\in\mathbb{V}_0^k$ can be written as $C=a_{-n}B$, and using \cref{recursiveY,normalord} we can always explicitly compute $Y(C,x)$.
Summarising what we have introduced so far, we have
\begin{itemize}
    \item a space of states $\mathbb{V}_0^k$,
    \item a vacuum vector $\vak\in \mathbb{V}_0^k$,
    \item a translation operator $T\in\End(\mathbb{V}_0^k)$,
    \item the state field map $Y(\cdot,x)$ as in \cref{ymap}.
\end{itemize}
This structure, together with some additional properties (see e.g. \cite[\textsection 2.4.4]{frenkelvertex}), defines a vertex algebra on $\mathbb{V}_0^k$.

\section{Construction of higher Hamiltonians}\label{Construction}
\subsection{The algebra of observables}
Let us introduce a set of complex numbers $\bm{k}=\{k_i\}_{i=1}^N$, where $N\in\ZZ_{>0}$ and $k_i\neq -2$ for all $i=1,\dots,N$. Consider the following tensor product of vacuum Verma modules
\begin{equation}
    \mathbb{V}_0^{\bm{k}} = \mathbb{V}_0^{k_1}\otimes\dots\otimes\mathbb{V}_0^{k_N}.
\end{equation}
This space can be interpreted as a module over the direct sum of $N$ copies of $\hf{sl}_2$.
Let us denote by $A^{(i)}\in\hf{sl}_2^{\oplus N}$ the copy of $A\in\hf{sl}_2$ in the $i$\textsuperscript{th} direct summand.

Let us denote by $\CC\vakk$ the one-dimensional vacuum representation of the ``positive modes'' Lie subalgebra $(\sl_2[t]\oplus\CC\mathsf{k})^{\oplus N}\subset \hf{sl}_2^{\oplus N}$, defined by $(\mathsf{k}^{(i)}-k_i)\vakk=0$ and $a_n^{(i)}\vakk=0$ for all $n\geq0$, $a\in\sl_2$ and $i=1,\dots,N$.
Therefore $\mathbb{V}_0^{\bm{k}}$ is the induced $\hf{sl}_2^{\oplus N}$-module, namely
\begin{equation}\label{VVM}
    \mathbb{V}_0^{\bm{k}} = U(\hf{sl}_2^{\oplus N})\otimes_{U(\sl_2[t]\oplus\CC\mathsf{k})^{\oplus N}}\CC\vakk.
\end{equation}

Repeating similar arguments to those of the previous sections, we can define $U_{\bm{k}}(\hf{sl}_2^{\oplus N})$ as the quotient of $U(\hf{sl}_2^{\oplus N})$ by the two-sided ideal generated by $\mathsf{k}^{(i)}-k_i$ for all $i=1,\dots,N$. We have the isomorphism
\begin{equation}
    U_{\bm{k}}(\hf{sl}_2^{\oplus N}) \cong U_{k_1}(\hf{sl}_2)\otimes U_{k_2}(\hf{sl}_2)\otimes\dots\otimes U_{k_N}(\hf{sl}_2).
\end{equation}
Thanks to this fact, $A^{(i)}\in\hf{sl}_2^{\oplus N}\subset U_{\bm{k}}(\hf{sl}_2^{\oplus N})$ can be presented as
\begin{equation}
    A^{(i)} = \mathds{1}\otimes \dots \otimes \mathds{1}\otimes A \otimes\mathds{1}\otimes \dots \otimes \mathds{1},
\end{equation}
where $A$ is acting as the $i$th tensor factor.
Again, we can introduce the left ideals $J_n^N\in U_{\bm{k}}(\hf{sl}_2^{\oplus N})$ generated by $a_r^{(i)}$ for all $r\geq n$, $a\in\mathfrak{sl}_2$ and $i=1,\dots,N$.
Let $\widetilde{U}_{\bm{k}}(\hf{sl}_2^{\oplus N})=\varprojlim U_{\bm{k}}(\hf{sl}_2^{\oplus N})/J_n^N$ be the inverse limit. This space is a complete topological algebra and
\begin{equation}
    \widetilde{U}_{\bm{k}}(\hf{sl}_2^{\oplus N}) \cong \widetilde{U}_{k_1}(\hf{sl}_2)\widehat{\otimes}  \cdots \widehat{\otimes} \widetilde{U}_{k_N}(\hf{sl}_2),
\end{equation}
where $\hat{\otimes}$ denotes the completed tensor product.
This space $ \widetilde{U}_{\bm{k}}(\hf{sl}_2^{\oplus N})$ is called the algebra of observables of the Gaudin model.

The tensor product $\mathbb V_0^{\bm k}$ is again a vertex algebra. The state-field map  $Y(\cdot,x):\mathbb{V}_0^{\bm{k}}\rightarrow \Hom(\mathbb{V}_0^{\bm{k}}, \mathbb{V}_0^{\bm{k}} ((x)))$ is just as in \cref{VAstructure} but decorated with the extra index $^{(i)}$.

Let us introduce the map $\Delta: \hf{sl}_2\hookrightarrow \hf{sl}_2^{\oplus N}$, which is the diagonal embedding of $\hf{sl}_2$ into $\hf{sl}_2^{\oplus N}$, defined as
\begin{equation}\label{diagaction}
    \Delta x = \sum_{i=1}^N x^{(i)}, \qquad \text{for all } x\in \hf{sl}_2.
\end{equation}
It extends to an embedding $\Delta:U(\hf{sl}_2)\hookrightarrow U(\hf{sl}_2^{\oplus N})\cong U(\hf{sl}_2)^{\otimes N}$. It is easy to check that
\begin{equation}
    [\Delta X_m, \Delta Y_n] = \Delta[X,Y]_{n+m}-n\delta_{n+m,0}(X|Y)\sum_{i=1}^N \mathsf{k}^{(i)},
\end{equation}
where $(\cdot|\cdot)$ is the usual Killing form as in \cref{killing}.

Therefore $\Delta$ descends to an embedding of the quotients $\Delta: U_{|\bm{k}|}(\hf{sl}_2)\hookrightarrow U_{\bm{k}}(\hf{sl}_2^{\oplus N})$, where $|\bm{k}|=\sum_{i=1}^N k_i$, and hence of their completions
\begin{equation}
    \Delta: \widetilde{U}_{|\bm{k}|}(\hf{sl}_2)\hookrightarrow \widetilde{U}_{\bm{k}}(\hf{sl}_2^{\oplus N}).
\end{equation}

\subsection{Invariant tensors on \texorpdfstring{$\mathfrak{sl}_2$}{TEXT}}
Let $\{I^a\}_{a=1}^3$ be a basis of $\sl_2$, and let $\{I_a\}_{a=1}^3$ be its dual basis with respect to the non-degenerate Killing form of \cref{killing}.
Let $f^{ab}{}_c$ denote the structure constants, so that
\begin{equation}
    [I^a,I^b]=f^{ab}{}_c I^c.
\end{equation}
Here and in what follows we employ the summation convention on Lie algebra indices.
Thanks to the non-degeneracy of the bilinear form, we may suppose our basis is chosen in such a way that
\begin{equation}\label{normalization}
    (I^a|I^b)=\delta^{ab}.
\end{equation}
By doing this, we no longer have to distinguish between upper and lower indices.
The structure constants are then
\begin{equation}\label{fabc}
    f^{ab}{}_c=f^{abc}=i\sqrt{2}\epsilon^{abc},
\end{equation}
where $\epsilon^{abc}$ is the usual Levi-Civita symbol.

(Concretely, in the defining representation we have
\begin{equation}
I^1 =\frac{1}{\sqrt{2}}\begin{pmatrix}0&1\\1&0\end{pmatrix}\qquad
I^2 =\frac{1}{\sqrt{2}}\begin{pmatrix}0&-i\\i&0\end{pmatrix}\qquad
I^3 =\frac{1}{\sqrt{2}}\begin{pmatrix}1&0\\0&-1\end{pmatrix}.
\end{equation}
Using \eqref{killing}, it is easy to check that \cref{normalization} holds.)

Recall that for any finite-dimensional Lie algebra $\mathring{\g}$, a tensor $t:\mathring{\g}\times\dots\times\mathring{\g}\rightarrow \CC$ is \emph{invariant} if
\begin{equation}\label{invariance1}
    t([a,x],y,\dots,z)+t(x,[a,y],\dots,z)+\dots+t(x,y,\dots,[a,z])=0, \qquad \text{for all } x\in \mathring{\g},
\end{equation}
or equivalently, if its components $t^{a_1\dots a_n} := t(I^{a_1} ,\dots,I^{a_n})$ satisfy
\begin{equation}\label{invariance2}
    f^{c a_1}{}_{b}t^{b a_2\dots a_n}
    +f^{c a_2}{}_{b}t^{a_1 b \dots a_n}
     + \dots + f^{c a_n}{}_{b}t^{a_1 a_2\dots b}
    =0,
\end{equation}
where the indices take values from $1$ to $\dim\mathring{\g}$.
In our case of $\sl_2$, the ring of invariant tensors is generated by $\delta^{ab}$ and $f^{a b c}$.
We shall need the following syzygy relations between them:
\begin{equation}
    \begin{gathered}
        f^{abc}f^{cde} = 2(\delta^{ae} \delta^{bd} -
 \delta^{ad} \delta^{be}), \qquad\qquad f^{abc}f^{abd} = -4\delta^{cd},\\
 f^{abc}\delta^{de} - f^{bcd}\delta^{ae} + f^{cda}\delta^{be} - f^{d a b}\delta^{ce} = 0.
    \end{gathered}
\end{equation}
Note in particular the last of these, which will play a crucial role in the explicit calculations of the next sections. It can also be generalized to higher rank tensors (see e.g. \cite[\textsection 369 F]{ito1993encyclopedic}).

\section{Quartic Hamiltonian}\label{quartic}
\subsection{Meromorphic states}\label{merstates}
Let us introduce a set $\{z_1,\dots,z_N\}$ of $N\in \ZZ_{>0}$ points $z_i\in \CC$ in the complex plane, chosen to be pairwise distinct, $z_i\neq z_j$ whenever $i\neq j$. For any element $A\in\hf{sl}_2$ we introduce the $\hf{sl}_2^{\oplus N}$-valued meromorphic functions
\begin{equation}\label{merfunc}
    A(z) := \sum_{i=1}^N\frac{A^{(i)}}{z-z_i}.
\end{equation}
We are allowed to take derivatives of such functions, which will be denoted by $A'(z)$ or, in general, for each $p\geq 0$,
\begin{equation}
    A^{\der p}(z) := \left(\frac{d}{dz}\right)^p A(z) = \sum_{i=1}^N (-1)^{p} p! \frac{A^{(i)}}{(z-z_i)^{p+1}}.
\end{equation}
Considering two of these functions with different spectral parameters, we get the following commutation relations
\begin{equation}\label{commpar}
\begin{split}
    [A^{\der p}(z),B^{\der q}(w)]=&(-1)^{p+1}(p+q)!\frac{[A,B](z)-[A,B](w)}{(z-w)^{p+q+1}}\\
    &+\sum_{k=1}^p(-1)^{p+1-k}\binom{p}{k}(p+q-k)!\frac{[A,B]^{\der k}(z)}{(z-w)^{p+q+1-k}}\\
    &-\sum_{k=1}^q(-1)^{p+1}\binom{q}{k}(p+q-k)!\frac{[A,B]^{\der k}(w)}{(z-w)^{p+q+1
    -k}}.
\end{split}
\end{equation}
By taking the limit $w\to z$, we get the commutation relations for the same spectral parameter, namely
\begin{equation}\label{comm1}
    [A^{\der p}(z),B^{\der q}(z)]=-[A,B]^{\der{p+q+1}}(z).
\end{equation}
We see that these $A^{\der p}(z)$, for $A\in \hf{sl}_2$ and $p\geq 0$, span a Lie algebra of $\hf{sl}_2^{\oplus N}$-valued meromorphic functions of $z$ with poles at the marked points.

It is helpful to be able to treat this as an abstract Lie algebra. Thus, let $\mathfrak L$ denote the Lie algebra over $\CC$ with basis consisting of $I_{n}^{a \der p}(z)$ and $\mathsf k^{\der p}(z)$, for $n\in \ZZ$, $p \in \ZZ_{\geq 0}$ and $a\in \{1,2,3\}$ with the non-vanishing Lie brackets given by
\begin{equation}\label{Irels}
    [I^{a \der p}_m(z),I^{b \der q}_n(z)] = - \frac{p!q!}{(p+q+1)!} (f^{ab}_c I^{c\der{p+q+1}}_{m+n}(z) - n \delta^{ab}\delta_{m+n,0}\mathsf k^{\der{p+q+1}}(z)).
\end{equation}
Let $\mathfrak L_+$ denote subalgebra generated by $I_{n}^{a \der p}(z)$ for $n\geq 0$, $p \in \ZZ_{\geq 0}$ and $a\in \{1,2,3\}$, and let
\begin{equation} \label{Vmodule}
\mathcal{V} := U(\mathfrak{L}) \otimes_{U(\mathfrak{L}_+)}\CC\vac
\end{equation}
denote the module over $\mathfrak L$ induced from the trivial one-dimensional module $\CC\vac$ over $\mathfrak L_+$.

We call the $\mathcal V$ the \emph{space of meromorphic states}. It is again a vertex algebra, with the state-field map as in \cref{VAstructure} but decorated with extra indices.
For each $z\in \CC\setminus\{z_1,\dots,z_N\}$, one has the homomorphism of Lie algebras $\mathfrak L \to \hf{sl}_2^{\oplus N}$ given by evaluating at $z$. It gives rise to a map $\mathcal V \to \mathbb V_0^{\bm k}$ of vertex algebras.

There is a bi-gradation of $\mathfrak{L}$ in which $X_{-n}^{\der p}(z)$ (for any $X\in \mathfrak{sl}_2$) has weight $(n,p+1)$ and $\mathsf k^{\der p}(z)$ has weight $(0,p+1)$. This yields a bi-gradation of $\mathcal{V}$
\begin{equation}\label{gradation}
    \mathcal{V} = \bigoplus_{n\geq 0, p\geq 0} \mathcal{V}_{n,p}.
\end{equation}
For each $n$, let $\mathcal{V}_n := \mathcal{V}_{n,n}$ denote the subspace of grade $(n,n)$. We call elements of $\mathcal{V}_n$ \emph{homogeneous meromorphic states of degree $n$}.

\subsection{Diagonal action of the zero modes of \texorpdfstring{$\hf{sl}_2$}{TEXT}}
There is an evident diagonal action of the Lie algebra $\hf{sl}_2$ on the $\mathfrak L$-module $\mathcal V$, defined in the same way as the action on $\mathbb V_0^{\bm k}$ in \cref{diagaction}.
In particular, for any $X\in\sl_2$, the zero modes stabilize each subspace $\mathcal{V}_{n,p}$, namely
\begin{equation}
    \Delta X_0 : \mathcal{V}_{n,p} \rightarrow \mathcal{V}_{n,p}.
\end{equation}

An important fact is that every state in $\mathcal{V}_n$ properly contracted with an $\sl_2$-invariant tensor vanishes under the diagonal action of the zero modes. This follows directly from the defining property of invariant tensors in \cref{invariance2}.
Let denote with $\mathcal{V}_n^{\sl_2}$ the invariant subspace. We can characterize this space for small $n$:

\begin{itemize}
    \item for $n=0$, $\mathcal{V}_0^{\sl_2}=\mathcal{V}_0=\mathbb C \vac$.
    \item for $n=1$, $\mathcal{V}_1^{\sl_2}= \{0\}$. Indeed, elements of $\mathcal V_1$ are of the form $t_a I^a_{-1}(z)\vac$. Such an element is in $\mathcal V^{\sl_2}_1$ if and only if $t_a$ are the components of an $\sl_2$-invariant tensor of rank 1. But there are no nonzero such tensors.
    \item for $n=2$, $\mathcal{V}_2^{\sl_2}$ has dimension 1 and it is spanned by the state
    \begin{equation}\label{S1}
        \varsigma_1(z) = \delta_{ab} I^{a}_{-1}(z)I^{b}_{-1}(z)\vac.
    \end{equation}
    \item for $n=3$, $\mathcal V^{\sl_2}_3$ has dimension 2 and it is spanned by the states
    \begin{equation}\label{3notworking}\begin{split}
        f^{abc}I^a_{-1}(z)I^b_{-1}(z)I^c_{-1}(z)\vac = &f^{abc}\frac{1}{2}(I^a_{-1}(z)I^b_{-1}(z) - I^b_{-1}(z)I^a_{-1}(z))I^c_{-1}(z)\vac\\
         =& f^{abc}f^{abd} I^{d\prime}_{-2}(z)I^c_{-1}(z) \vac\\
         = & -4 I^{c\prime}_{-2}(z)I^c_{-1}(z) \vac,
    \end{split}
    \end{equation}
    and
        \begin{equation}I^{c}_{-2}(z)I^c_{-1}(z)\mathsf{k}(z)\vac. \end{equation}
    \item for $n=4$, $\mathcal{V}_4^{\sl_2}$ has dimension 14. Below, we will make use of the following explicit choice of basis:
    \begin{equation}\label{genquar}
\begin{gathered}
    \mathsf{v}_1 := \delta^{(ab}\delta^{cd)}I^{a}_{-1}(z)I^{b}_{-1}(z)I^{c}_{-1}(z)I^{d}_{-1}(z), \qquad  \mathsf{v}_2 := f^{abc}I^{a}_{-2}(z)I^{b\prime}_{-1}(z)I^{c}_{-1}(z),\\
    \mathsf{v}_3 := I^{a\prime\prime}_{-3}(z)I^{a}_{-1}(z),\qquad \mathsf{v}_4:=I^{a}_{-3}(z)I^{a\prime\prime}_{-1}(z), \qquad \mathsf{v}_5:=I^{a\prime\prime}_{-2}(z)I^{a}_{-2}(z),\\ \mathsf{v}_6:=I^{a\prime}_{-3}(z)I^{a\prime}_{-1}(z),\qquad \mathsf{v}_7:=I^{a\prime}_{-2}(z)I^{a\prime}_{-2}(z),\qquad
    \mathsf{v}_8:=I^{a}_{-3}(z)I^{a}_{-1}(z)\mathsf{k}^{\prime}(z),\\ \mathsf{v}_9:=I^{a}_{-2}(z)I^{a}_{-2}(z)\mathsf{k}^{\prime}(z) , \qquad \mathsf{v}_{10}:=I^{a\prime}_{-3}(z)I^{a}_{-1}(z)\mathsf{k}(z),\\ \mathsf{v}_{11}:=I^{a}_{-3}(z)I^{a\prime}_{-1}(z)\mathsf{k}(z),\qquad
    \mathsf{v}_{12}:=I^{a\prime}_{-2}(z)I^{a}_{-2}(z)\mathsf{k}(z),\\ \mathsf{v}_{13}:=I^{a}_{-3}(z)I^{a}_{-1}(z)\mathsf{k}(z)^2,\qquad \mathsf{v}_{14}:=I^{a}_{-2}(z)I^{a}_{-2}(z)\mathsf{k}(z)^2.
\end{gathered}
\end{equation}
\end{itemize}
Note that to write these terms we have to choose an ordering prescription. In this work we sort level first in ascending order from left to right and after that, for a given level, we sort derivatives in descending order from left to right. For example $f^{abc}I^{a}_{-2}I^{b\prime\prime}_{-2}I^{c\prime}_{-3} = f^{abc}I^{c\prime}_{-3}I^{b\prime\prime}_{-2}I^{a}_{-2} + \text{terms obtained from commutations}$.

\subsection{Top terms}
We can see from the above construction that in the case $n=2$ and $n=4$, there is a particular state, that we will call \emph{top term}, which is the state in $\mathcal{V}_n^{\sl_2}$ that contains exactly $n$ generators:
\begin{equation}
    \delta_{ab} I^{a}_{-1}(z)I^{b}_{-1}(z)\vac, \qquad
    \delta_{(ab}\delta_{cd)}I^{a}_{-1}(z)I^{b}_{-1}(z)I^{c}_{-1}(z)I^{d}_{-1}(z) \vac.
\end{equation}
We do not have such state for $n=3$, because we can always use the commutation relations to reduce the number of generators, as shown in \cref{3notworking}.
This pattern continues, as we now describe.

Notice that the universal enveloping algebra $U(\mathfrak L)$ has an increasing filtration
\begin{equation}
    F_0U(\mathfrak{L}) \subseteq F_1U(\mathfrak{L}) \subseteq \dots \subseteq U(\mathfrak{L}),
\end{equation}
in which the generators $I^{a[p]}_n(z)$ count as $+1$ and the generators $\mathsf k^{[p]}(z)$ count as $0$, cf. the commutation relations of $\mathfrak L$ in \cref{Irels}.
For example $I^a_{-1}(z)I^{a \prime}_{-2}(z)\in F_2$, and $I^a_{-1}(z)I^{a \prime}_{-2}(z)\mathsf{k}(z)\in F_2$ as well.
It gives rise to a corresponding filtration, $F_0\mathcal V \subseteq F_1\mathcal V \subseteq \dots \subseteq \mathcal V$, on the space $\mathcal V$ of meromorphic states.

Observe that if $\mathsf{v}\in \mathcal V_N$ then $\mathsf{v}\in F_N \V_N$.  We see that
\begin{align}
\mathsf{v} &\equiv t_{i_1\dots i_N}I^{i_1}_{-1}(z)\dots I^{i_N}_{-1}(z)\vac \mod F_{N-1}\V_N\nonumber\\
  &\equiv t_{(i_1\dots i_N) }I^{i_1}_{-1}(z)\dots I^{i_N}_{-1}(z)\vac \mod F_{N-1}\V_N,
\end{align}
for some $\sl_2$ tensor $t_{i_1,\dots,i_N}$, where the brackets around the indices denote the operation of symmetrization,
\begin{equation}
    t_{(i_1,\dots,i_n)}=\frac{1}{n!}\sum_{\sigma\in\mathcal{S}_n}t_{\sigma(i_1)\dots\sigma(i_n)}
\end{equation}
(and we may symmetrize without loss of generality because the non-symmetric pieces fall into $F_{N-1}$, as for example in \cref{3notworking}).
Let us call $t_{(i_1\dots i_N) }I^{i_1}_{-1}(z)\dots I^{i_N}_{-1}(z)\vac$ the \emph{top term} of the state $\mathsf{v}\in \mathcal V_N$.

If this state $\mathsf v\in \V$ is $\sl_2$-invariant, $\mathsf v\in \V^{\sl_2}$, then $t_{(i_1\dots i_N) }$ is a symmetric invariant tensor. Nonzero such tensors exist only in even degrees, and up to rescaling they are, explicitly,
\begin{equation}\label{invtens}
\begin{gathered}
    t_{i_1,i_2} = \delta_{i_1i_2}\\
    t_{i_1,i_2,i_3,i_4} = \delta_{(i_1i_2}\delta_{i_3i_4)}\\
    t_{i_1,i_2,i_3,i_4,i_5,i_6} = \delta_{(i_1i_2}\delta_{i_3i_4}\delta_{i_5i_6)}\\\nonumber
    \dots
\end{gathered}
\end{equation}

In what follows, our interest is in meromorphic states $\mathsf v\in \V^{\sl_2}$ that have nonzero top term (in other words states whose principal symbol has maximal degree) and that are $\sl_2$-invariant.

\subsection{Singular vectors up to twisted derivative}\label{twistedder}
Let us define the twisted derivative operator of degree $j\in\ZZ$ with respect to the spectral parameter $z$,
\begin{equation}\label{TD}
    D_z^{(j)} = \left( \partial_z-\frac{j}{2}\mathsf{k}(z)\right).
\end{equation}
Note that this operator sends $\V_{n,p} \to \V_{n,p+1}$ in the bigradation we introduced above.

We will say that a vector $\mathsf{v}\in \mathcal{V}^{\sl_2}_{n}$ is \emph{singular up to twisted derivatives} if for all $x\in\sl_2$  we have
\begin{equation}\label{singularity}
    \Delta x_m \mathsf{v} = 0 \qquad \textup{mod}\: D^{(n-1)}_z\mathcal{V}_{n-m,n-1}.
\end{equation}
for \emph{all} non-negative modes $x_m$, $m\geq 0$. This defines a subspace
\begin{equation}
    \V_n^{\text{sing}} \subset \V_n^{\sl_2}
\end{equation}
of vectors singular up to twisted derivatives.
\begin{proposition}\label{singS1}
The space of singular vectors $\mathcal{V}^{\text{sing}}_2$ is spanned by the quadratic state $\varsigma_1(z)$ defined in \eqref{S1}.
\end{proposition}

\hrule
\begin{proof}
We need to show that
\begin{equation}
    \Delta I^r_{k} \varsigma_{1}(z) = 0 \quad \mod D^{(1)}_z \mathsf{G}^r_{k}(z),
\end{equation}
for some meromorphic states $\mathsf{G}_k^{r}\in\mathcal{V}_{2-k,1}$, for all $k\geq 0$ and $r=1,2,3$.
For $k=0$ there is nothing to check since $\Delta I^r_{0} \varsigma_{1}(z) = 0$ identically, by the definition of $\mathcal{V}^{\sl_2}_2$.
It is enough to check the action of the first modes $I^r_1$, since any higher modes can be expressed in terms of their brackets, \textit{i.e.} $I_{2}^r= -\frac{1}{4} f^{rbc} [I^b_1,I^c_1]$ etc.
From direct calculations we get that
\begin{equation}
    \Delta I^r_{1} \varsigma_{1}(z) = D^{(1)}_z \mathsf{G}^r_{1}(z),
\end{equation}
where
\begin{equation}\label{G11}
    \mathsf{G}^r_{1}(z) = -4 I^r_{-1}(z)\vac.
  \end{equation}
  \hrule
\end{proof}

More non-trivially, for $n=4$ we have the following result.
\begin{proposition}\label{singS3}
The space of singular vectors $\mathcal{V}^{\text{sing}}_4$ is of dimension 7. A choice of basis is given by the state
\begin{equation}\label{S3}
        \begin{split}
    \varsigma_3(z) =\Big[&\delta_{(ab}\delta_{cd)}I^{a}_{-1}(z)I^{b}_{-1}(z)I^{c}_{-1}(z)I^{d}_{-1}(z)+\frac{20}{3}  f_{abc}I^{a}_{-2}(z)I^{b\prime}_{-1}(z)I^{c}_{-1}(z)\\
    &+\frac{40}{9}I^{a}_{-3}(z)I^{a\prime\prime}_{-1}(z)-\frac{20}{3}I^{a\prime\prime}_{-2}(z)I^{a}_{-2}(z)
    +\frac{40}{9}I^{a\prime}_{-3}(z)I^{a\prime}_{-1}(z)\\
    &\phantom+ \qquad -\frac{10}{3}I^{a\prime}_{-2}(z)I^{a\prime}_{-2}(z)-\frac{20}{3}I^{a}_{-3}(z)I^{a}_{-1}(z)\mathsf{k}^{\prime}(z)\Big]\vac,
    \end{split}
\end{equation}
together with the double translate state
\begin{equation}\label{vectrans}
T^2\Big(I^{a\prime\prime}_{-1}(z)I^{a}_{-1}(z)\vac-I^{a\prime}_{-1}(z)I^{a\prime}_{-1}(z)\vac-\frac{3}{4}I^{a}_{-1}(z)I^{a}_{-1}(z)\mathsf{k}^{\prime}(z)\vac\Big)
\end{equation}
and the following twisted derivative states
\begin{gather}\label{vecTD}
    D^{(3)}_z \Big( I^{a}_{-3}(z)I^{a\prime}_{-1}(z)\vac\Big),\quad
    D^{(3)}_z \Big( I^{a\prime}_{-3}(z)I^{a}_{-1}(z)\vac\Big),\quad
    D^{(3)}_z \Big( I^{a\prime}_{-2}(z)I^{a}_{-2}(z)\vac\Big), \\
    D^{(3)}_z \Big( I^{a}_{-3}(z)I^{a}_{-1}(z)\mathsf{k}(z)\vac\Big), \quad
    D^{(3)}_z \Big( I^{a}_{-2}(z)I^{a}_{-2}(z)\mathsf{k}(z)\vac\Big).\nonumber
\end{gather}
\end{proposition}

\begin{mproof}
Let  $s(z)\in \V_4^{\sl_2}$. We may write it in our basis \eqref{genquar},
\begin{equation}
    s(z) = \sum_{i=1}^{14} \xi_i \mathsf v_i(z)
\end{equation}
for some coefficients $\xi_i\in\CC$ with $i=1,\dots,14$, and ask what conditions the requirement of being singular up to twisted derivatives, \eqref{singularity}, places on these coefficients. It is enough to demand that
\begin{equation}
    \Delta I^r_k s(z) =0 \mod  D^{(3)}_z \mathsf{G}^r_{k}(z)
\end{equation}
for some meromorphic states $\mathsf{G}^r_{k}(z) \in \V_{4-k,3}$, for all $k\geq 0$ and $r=1,2,3$. For zero modes there is nothing to check since $\Delta I^a_0 s(z) =0$ exactly, by definition of $\V_4^{\sl_2}$. It is then enough to check the action of first modes, $I^r_1$, since any higher modes can be expressed in terms of their brackets, $I^r_2= -\frac{1}{4} f^{rbc} [I^b_1,I^c_1]$ etc. So we are to check under what conditions
\begin{equation}
    \Delta I^r_1 s(z) = D^{(3)}_z \mathsf{G}^r_{1}(z)
\end{equation}
for some $\mathsf{G}^a_{1}(z) \in \V_{3,3}$. By direct calculation, one finds that solutions exist precisely if the coefficient $\xi_i$ obey the relations
\begin{equation}\label{coeffs31}
    \begin{gathered}
    \xi_2= \frac{20}{3}\xi_1, \qquad
    \xi_3 = \frac{20}{3}\xi_1 -\xi_4 +2\xi_5+\xi_6-2\xi_7, \qquad
    \xi_9 = -\frac{5}{4}\xi_1-\frac{3}{8}\xi_5+\frac{3}{8}\xi_7-\frac{2}{3}\xi_{14}, \qquad\\
    \xi_{10} = \frac{5}{3}\xi_1+\frac{3}{2}\xi_4-\frac{3}{2}\xi_5-\frac{3}{2}\xi_6+\frac{3}{2}\xi_7+\xi_8,\qquad
    \xi_{11} = \frac{55}{3}\xi_1-\frac{3}{2}\xi_4+\frac{3}{2}\xi_5-\frac{3}{2}\xi_7+\xi_8,\qquad \\
    \xi_{12} = -\frac{15}{2}\xi_1-\frac{3}{4}\xi_5-\frac{3}{4}\xi_7-\frac{4}{3}\xi_{14},\qquad
    \xi_{13}=-\frac{55}{4}\xi_1-\frac{9}{8}\xi_5+\frac{9}{8}\xi_7-\frac{3}{2}\xi_8.
    \end{gathered}
\end{equation}
When they do obey these relations, the required functions $\mathsf{G}^r_{1}(z)$ are given by
\begin{equation}\label{G13}
\begin{split}
    \mathsf{G}^r_{1}(z)  = \Big[ & \rho_1  I^{(a}_{-1}(z)I^r_{-1}(z)I^{a)}_{-1}(z)\\
    &+ f^{r a b}\left(\rho_2 I^a_{-2}(z)I^{b\prime}_{-1}(z) + \rho_3 I^{a\prime}_{-2}(z)I^{b}_{-1}(z)+ \rho_4 I^a_{-2}(z)I^{b}_{-1}(z)\mathsf{k}(z)\right)\\
    &+ \rho_5 I^{r}_{-3}(z)\mathsf{k}(z)^2 + \rho_6 I^{r}_{-3}(z)\mathsf{k}^{\prime}(z) + \rho_7 I^{r\prime}_{-3}(z)\mathsf{k}(z) + \rho_8 I^{r\prime\prime}_{-3}(z)\Big] \vac,
\end{split}
\end{equation}
where
\begin{equation*}
\begin{gathered}
    \rho_1 = -\frac{8}{3}\xi_1, \qquad \rho_2 =\frac{20}{3}\xi_1 -\xi_4 +  \xi_5,
    \qquad \rho_3 = \xi_4-\xi_5-\xi_6+2\xi_7,\\
    \rho_4 = -\frac{55}{6}\xi_1 -\frac{3}{4}\xi_5+\frac{3}{4}\xi_7-\xi_8-\frac{4}{3}\xi_{14} \qquad  \rho_5 = \frac{55}{6}\xi_1+\frac{3}{4}\xi_5-\frac{3}{4}\xi_7+\xi_8,
    \qquad \rho_6 = \xi_4,\\
    \rho_7 = 5\xi_1 -\xi_4+\frac{3}{2}\xi_5+\xi_6-\frac{3}{2}\xi_7+\frac{8}{3}\xi_{14}, \qquad \rho_8 =-\frac{100}{9}\xi_1-\frac{4}{3}\xi_5-\frac{2}{3}\xi_7.
\end{gathered}
\end{equation*}
The basis reported in the proposition can be obtained by the one defined by the restrictions \eqref{coeffs31} by a change of basis.
\end{mproof}

The proposition above is in agreement with the calculation of the quartic Hamiltonian density $S_4(z)$ (the analogue of our $\varsigma_3(z)$), recently presented in \cite{LKT22}. In the present conventions, the latter is given by
\begin{equation}\begin{split}
    S_4(z) = \Big[& \delta^{(ab}\delta^{cd)} I^{a}_{-1}(z)I^{b}_{-1}(z)I^{c}_{-1}(z)I^{d}_{-1}(z) +\frac{20}{3}f^{abc} I^{a}_{-2}(z)I^{b\prime}_{-1}(z)I^{c}_{-1}(z) \\
    &-\frac{40}{9}I^{a\prime\prime}_{-3}(z)I^{a}_{-1}(z) -\frac{140}{9} I^{a\prime\prime}_{-2} (z) I^{a}_{-2} (z) +\frac{40}{3}I^{a\prime}_{-3} (z) I^{a\prime}_{-1} (z) \\
    &\phantom+\qquad-\frac{10}{3}I^{a\prime}_{-2} (z) I^{a\prime}_{-2} (z) + 5 I^{a}_{-2} (z) I^{a}_{-2} (z) \mathsf{k}(z)^2\Big]\vac
\end{split}
\end{equation}
and it does\footnote{To match conventions, note that for us
\begin{equation}
    \delta_{(ab}\delta_{cd)} = \frac{1}{3} \left(
    \delta_{ab}\delta_{cd} + \delta_{ac}\delta_{bd} + \delta_{ad}\delta_{bc}\right)
\end{equation}
and in \cite{LKT22} the tensor called $\tau_3^{abcd}$ is given by
\begin{equation}
    \tau_3^{abcd} = \frac 1{16} \left(\delta_{ab}\delta_{cd} + \delta_{ac}\delta_{bd} + \delta_{ad}\delta_{bc}\right).
\end{equation}
We thank Sylvain Lacroix for clarifying discussions on this point.} indeed lie in the space $\mathcal{V}^{\text{sing}}_4$.

\subsection{Hamiltonian densities}

Now, to state the main result of the paper, we need two reintroduce rational functions of two different spectral parameters, $z$ and $w$, cf. \cref{commpar} and \cref{comm1}.

Recall the Lie algebra $\mathfrak L\equiv \mathfrak L_{(z)}$ over $\CC$ from \cref{merstates}. Let $\mathfrak L_{(z,w)}$ be the Lie algebra with generators $I_n^{a[p]}(z)$, $I_n^{a[p]}(w)$, $\mathsf k^{[p]}(z)$ and $\mathsf k^{[p]}(w)$ for $a=1,2,3$, $n\in \mathbb Z$ and $p\in \mathbb Z_{\geq 0}$ and commutation relations
\begin{equation}
\begin{split}
    [I_m^{a\der p}(z),I_n^{b\der q}(w)]= & (-1)^{p+1}(p+q)!\frac{f^{ab}_c(I^{c}_{m+n}(z)-I^c_{m+n}(w)) - n \delta_{m+n,0}\delta^{ab}(\mathsf{k}(z) - \mathsf{k}(w) )}{(z-w)^{p+q+1}}\\
    &+\sum_{j=1}^p(-1)^{p+1-j}\binom{p}{j}(p+q-j)!\frac{f^{ab}_cI^{c \der j}_{m+n}(z)-n\delta_{n+m,0}\delta^{ab}\mathsf k^{\der j}(z)}{(z-w)^{p+q+1-k}}\\
    &-\sum_{k=1}^q(-1)^{p+1}\binom{q}{k}(p+q-k)!\frac{f^{ab}_cI^{c \der j}_{m+n}(w)-n\delta_{n+m,0}\delta^{ab}\mathsf k^{\der j}(w)}{(z-w)^{p+q+1
    -k}},
\end{split}
\end{equation}
together with \cref{Irels} for the generators with parameter $z$ and the analogue with parameter $w$.
This Lie algebra $\mathfrak L_{(z,w)}$ and its modules are defined over the ground ring $\CC[(z-w)^{-1}]$ of polynomials in $(z-w)^{-1}$.
We have the vertex algebra $\mathcal V_{(z,w)}$ defined analogously to \cref{Vmodule} and the two obvious embedding maps of vertex algebras $\mathcal V \hookrightarrow \mathcal V_{(z,w)}$, which we write as $\mathsf v \mapsto \mathsf v(z)$ and $\mathsf v \mapsto \mathsf v (w)$.

Moreover, there is a natural notion of ``expanding around $z=w$''. Namely, there is a homomorphism $\mathfrak L_{(z,w)} \to \mathfrak L_{(z)}((w-z))$ of Lie algebras over $\CC[(z-w)^{-1}]$ defined by
\begin{equation}\label{taylor}
    I^{a[p]}_m(w) = I^{a \der p}_m (z) + I^{a \der{p+1}}_m(z) (w-z) + \frac{1}{2}I^{a\der{p+2}}_m(z)(w-z)^2 +\dots
\end{equation}
which is motivated by considering the Taylor expansion $\iota_{w-z} A(w)$ of the function $A(w)$ from \cref{merfunc}.
This gives rise to a map $\mathcal V_{(z,w)} \to \mathcal V_{(z)}((w-z))$.
We say a state $\mathsf v\in \mathcal V_{(z,w)}$ is \emph{regular at $z=w$ modulo translates} if there exists $Z\in\mathcal V_{(z,w)}$ such that the image of $\mathsf v- TZ$ under this map has no singularities in $(z-w)$.

Recall from \cref{S1,S3} the definitions of the quadratic state $\varsigma_1\in\mathcal{V}_{2}^{\textup{sing}}$ and of the vector  $\varsigma_3\in\mathcal{V}_{4}^{\textup{sing}}$, respectively.
\begin{theorem}\label{th2}
The elements $\varsigma_1\in\mathcal{V}_{2}^{\textup{sing}}$ and $\varsigma_3\in\mathcal{V}_{4}^{\textup{sing}}$ obey the relations
\begin{subequations}
    \begin{align}
        \varsigma_1(z)_{(0)}\varsigma_1(w) &= (D_z^{(1)}-D_w^{(1)})\mathsf{A}_{1,1}(z,w) + T \mathsf{B}_{1,1}(z,w),\label{thS1S1} \\
        \varsigma_1(z)_{(0)}\varsigma_3(w) &= (3D_z^{(1)}-D_w^{(3)})\mathsf{A}_{1,3}(z,w) + T \mathsf{B}_{1,3}(z,w),\label{thS1S3}\\
        \varsigma_3(z)_{(0)}\varsigma_1(w) &= (D_z^{(3)}-3D_w^{(1)})\mathsf{A}_{3,1}(z,w) + T \mathsf{B}_{3,1}(z,w),\label{thS3S1}\\
        \varsigma_3(z)_{(0)}\varsigma_3(w) &= (3D_z^{(1)}-3D_w^{(3)})\mathsf{A}_{3,3}(z,w) + T \mathsf{B}_{3,3}(z,w),\label{thS3S3}
    \end{align}
\end{subequations}
where $\mathsf{A}_{i,j}(z,w)$ and $\mathsf{B}_{i,j}(z,w)$ are elements of $\mathcal V _{(z,w)}$.
Moreover, $\mathsf{A}_{i,j}(z,w)$, $i,j\in\{1,3\}$, are regular at $z=w$ modulo translates.
\end{theorem}

\begin{mproof}
The two statements of the theorem follow from direct calculations.
In particular, when $m=n=1$, we get
\begin{equation}\label{resultS1S1}
    \mathsf{A}_{1,1}(z,w) = \frac{8}{z-w} I^{a}_{-2}(z)I^{a}_{-1}(w)\vac, \qquad \mathsf{B}_{1,1}(z,w) = \frac{8}{(z-w)^2} I^{a}_{-1}(z)I^{a}_{-1}(w)\vac .
\end{equation}
We have computed $\mathsf{A}_{1,3}(z,w)$, $\mathsf{B}_{1,3}(z,w)$, $\mathsf{A}_{3,3}(z,w)$ and $\mathsf{B}_{3,3}(z,w)$ explicitly, with the aid of the computer algebra system FORM \cite{FORM,FORM1}. The expressions for $\mathsf{A}_{1,3}(z,w)$ and $\mathsf{B}_{1,3}(z,w)$, are given in \cref{AppA}. The expressions for $\mathsf{A}_{3,3}(z,w)$ and $\mathsf{B}_{3,3}(z,w)$ are extremely lengthy (more that 500 terms in total), and we do not reproduce them here.

Once the expression of $\varsigma_1(z)_{(0)}\varsigma_3(w)$ is known, \textit{i.e.} the functions  $\mathsf{A}_{13}(z,w)$ and $\mathsf{B}_{13}(z,w)$ are found, it can be shown that the theorem is automatically satisfied for the product $\varsigma_3(z)_{(0)}\varsigma_1(w)$.
This comes from the property of the $n$th product between two states $a,b$ of a vertex algebra, namely
\begin{equation}\label{skewsymm}
    a_{(n)}b =- \sum_{k=0}^{\infty}\frac{1}{k!}(-1)^{k+n}T^{k}(b_{(n+k)}a).
\end{equation}
Therefore, by swapping two states in a 0th product, we obtain  $\mathsf{A}_{31}(z,w)=-\mathsf{A}_{13}(w,z)$ and a series of terms which are nothing but translates and therefore can be absorbed in the definition of $\mathsf{B}_{31}(z,w) = -\mathsf{B}_{13}(w,z) + \sum_{k=0}^{\infty}(-1)^{k}T^{k}(\varsigma_1(w)_{(k+1)}\varsigma_3(z))$.

To prove the second part of the theorem one can expand according to \cref{taylor} and the result follows from direct calculation.
\end{mproof}

Having established this statement for the particular choice of quartic density $\varsigma_3$, we automatically get the following property for \emph{any} element of $\mathcal{V}^{\text{sing}}_4$. It is a slightly weaker property, because the condition on the twisted derivative terms on the right hand side is less rigid. As we shall see in \cref{commutingham} below, it is sufficient for defining consistent Hamiltonians.
\begin{corollary}\label{cor}
For \emph{any} element $\mathsf v_3 \in \mathcal{V}^{\text{sing}}_4$, one has
  \begin{align}
        \varsigma_1(z)_{(0)}\mathsf v_3(w) &= D_z^{(1)} \mathsf A^{\rom{1}}_{1,3}(z,w) + D_w^{(3)}\mathsf{A}^{\rom{2}}_{1,3}(z,w) + T \mathsf{B}_{1,3}(z,w),\\
        \mathsf v_3(z)_{(0)}\varsigma_1(w) &= D_z^{(3)} \mathsf A^{\rom{1}}_{3,1}(z,w) + D_w^{(1)}\mathsf{A}^{\rom{2}}_{3,1}(z,w) + T \mathsf{B}_{3,1}(z,w),\\
        \mathsf v_3(z)_{(0)}\mathsf v_3(w) &= D_z^{(3)} \mathsf A^{\rom{1}}_{3,3}(z,w) + D_w^{(3)}\mathsf{A}^{\rom{2}}_{3,3}(z,w) + T \mathsf{B}_{3,3}(z,w).
    \end{align}
where $\mathsf{A}^{\rom{1},\rom{2}}_{i,j}(z,w)$ and $\mathsf{B}_{i,j}(z,w)$ are elements of $\mathcal V _{(z,w)}$.
Moreover, $\mathsf{A}^{\rom{1},\rom{2}}_{ij}(z,w)$, $i,j\in\{1,3\}$, are regular at $z=w$ modulo translates.
\end{corollary}

\begin{mproof}
We already know from \cref{th2} that there exists an element, $\varsigma_3(z)$, satisfying these relations.
But we saw in \cref{singS3} that every element $\mathsf v_3(z)$ of $\mathcal V^{\text{sing}}_4$ is proportional to $\varsigma_3(z)$ up to the addition of certain translates and twisted derivatives.

It follows from the property \eqref{skewsymm} that if we add to $\varsigma_3(z)$ any translate then the statement of \cref{th2} still holds, the only difference being a re-definition of the states $\mathsf{B}(z,w)$. And it is evident that, if we add to $\varsigma_3(z)$ any linear combination of the twisted derivatives in \cref{vecTD} then the resulting vector $\mathsf v_3(z)$ still obeys the weaker relations given above. (One might worry about introducing singularities at $z=w$, but note that for \emph{any} meromorphic states $a(z)$ and $b(z)$, the product $a(z)_{(0)}b(w)$ is regular at $z=w$, as is manifest if we expand $b(w)$ about $w=z$ in the spectral plane before taking the vertex-algebra product:
$a(z)_{(0)} b(w) = a(z)_{(0)} \left(b(z) + (w-z) b'(z) + \dots \right) = a(z)_{(0)} b(z) + (w-z) a(z)_{(0)} b'(z) + \dots$.)
\end{mproof}

\subsection{Gaudin Hamiltonian}
Let us define the following state at non-critical level, \textit{i.e.} $k_i\neq -2$,
\begin{equation}
s_1(z) = \frac{1}{2} \Big( \varsigma_1(z) + 4 D^{(1)}_z\omega(z) \Big) \in \VV_0^{\bm k},
\end{equation}
where $\varsigma_1(z)$ is now the image in $\mathbb V_0^{\bm k}$ of the density defined in \eqref{S1} and where
\begin{equation}
    \omega(z) :=\sum_{i=1}^N \frac{1}{z-z_i}\left(\frac{1}{2(k_i+2)}I^{a (i)}_{-1}I^{a (i)}_{-1}\vakk\right),
\end{equation}
the term in the brackets being the Segal-Sugawara vector at site $i$.

It is possible to show (see \cite{LVYcubic}), that the operator $(s_1(z))_{(0)}$ is the image in $\tilde{U}_{\bm k} (\hf{sl}_2^{\oplus N})$ of
\begin{equation}
    \sum_{i=1}^N\frac{\mathscr{C}^{(i)}}{(z-z_i)^2}+\sum_{i=1}^N\frac{\mathscr{H}_i}{z-z_i}\in \tilde U (\g^{\oplus N}),
\end{equation}
where
\begin{equation}
    \mathscr{C}^{(i)} :=(\mathsf k^{(i)}+2)\mathsf{d}^{(i)} + \frac{1}{2}I^{a (i)}_{0}I^{a (i)}_{0}+\sum_{n> 0 }I^{a (i)}_{-n}I^{a (i)}_{-n}
    \end{equation}
is the $i$th copy of the quadratic Casimir operator of $\g$ in $\tilde{U}(\g^{\oplus N})$ and $\mathscr{H}_i$ are the Hamiltonians in \eqref{hamiltonian}.
\begin{theorem}
Given the images in $\VV_0^{\bm k}$ of the densities $\varsigma_i$, $i\in\{1,3\}$, we have
\begin{equation}
    s_1(z)_{(0)}\varsigma_i(w) = -\frac{1}{2}D^{(1)}_w\mathsf{A}_{1,i}(z,w)+T\Big( \frac{1}{2}\mathsf{B}_{1,i}(z,w)+2 D^{(1)}_z \frac{\varsigma_i(w)}{z-w}\Big),
\end{equation}
with $\mathsf{A}_{1,i}(z,w)$ and $\mathsf{B}_{1,i}(z,w)$ being the images in $\VV_0^{\bm k}$ of the meromorphic states in \cref{th2}.
\end{theorem}

\begin{mproof}
The result follows from direct calculations, using the definitions of $\mathsf{A}_{1,1}$, $\mathsf{B}_{1,1}$, $\mathsf{A}_{1,3}$, $\mathsf{B}_{1,3}$ in \eqref{resultS1S1} and \cref{AppA}, respectively.
\end{mproof}

As we will see in the next section, this requirement is sufficient to ensure the commutativity of local Hamiltonians, arising from the densities $\varsigma_1(z)$ and $\varsigma_3(z)$, with the usual quadratic Gaudin Hamiltonians which define the model.

\subsection{Commuting Hamiltonians}\label{commutingham}
In this section, we will simply recall the ideas presented in \cite{LVYaffine}.
Consider two states $X,Y\in \mathbb{V}_0^{\bm{k}}$ and their formal zero modes $X_{(0)}, Y_{(0)}\in\widetilde{U}_{\bm{k}}(\hf{sl}_2^{\oplus N})$. It is possible to show the following vertex algebra identity (see e.g. \cite{frenkelvertex})
    \begin{equation}\label{commzero}
        [X_{(0)}, Y_{(0)}]=(X_{(0)} Y)_{(0)}.
    \end{equation}
This means that if one is able to find a family of operators whose 0th product vanishes (or that can be expressed as a translation, since $(TZ)_{(0)}=0$ by definition), then their formal zero modes form a commutative subalgebra of the algebra of observables $\widetilde{U}_{\bm{k}}(\widehat{\mathfrak{sl}}^{\oplus N}_2)$.

The meromorphic function which is obtained by acting with $\mathsf{k}(z)$ on the module $\VV_0^{\bm k}$, \emph{i.e.} setting the central elements to numbers,
\begin{equation}\label{twistfunc}
    k (z)=\sum_{i=1}^N \frac{k_i}{z-z_i},
\end{equation}
has a special role and it is called the \emph{twist function} of the model.
Let us define also
    \begin{equation}\label{hyp}
        \mathscr{P}(z)\defeq\prod_{j=1}^{N}(z-z_j)^{k_j}.
    \end{equation}
The function $\mathscr P^{1/2}$ is multi-valued. It becomes single-valued on a certain multi-sheeted cover of $\CC\setminus\{z_1,\dots,z_N\}$. Let $\gamma$ be any closed contour in this cover. For example, $\gamma$ could be the lift to this cover of a Pochhammer contour in $\CC\setminus\{z_1,\dots,z_N\}$ around any two of the marked points.
Then $\mathscr P^{n/2}$, for any integer $n$, is single-valued along $\gamma$, and one can introduce the integral
    \begin{equation}
         \int_{\gamma} \mathscr{P}(z)^{-n/2} f(z)\diff z
    \end{equation}
which has the fundamental property that, for any meromorphic function $f(z)$ which is non-singular along $\gamma$
    \begin{equation}\label{Pfunction}
        \int_{\gamma} \mathscr{P}(z)^{-n/2}\left( \frac{d}{dz} - \frac{n}{2} k(z) \right) f(z)\diff z= \int_{\gamma} \frac{d}{dz} (\mathscr{P}^{-n/2}f(z)) \diff z=0,
    \end{equation}
cf. \cref{TD}. (See e.g. \cite{LVYaffine} for the details.)

Let us now define the following object in $\widetilde{U}_{\bm{k}}(\hf{sl}_2^{\oplus N})$,
\begin{equation}
    Q_n^{\gamma} = \int_{\gamma} \mathscr{P}(z)^{-n/2} \varsigma_n(z)_{(0)} \diff z ,
\end{equation}
for $n=1,3$,
where $\varsigma_n(z)$ are now the images in $\VV_0^{\bm k}$ of the densities we have defined in the previous section.
\begin{proposition}\label{prop}
The operators $Q_n^{\gamma}\in\widetilde{U}_{\bm{k}}(\hf{sl}_2^{\oplus N})$ commute amongst themselves, with the generators of $\hf{sl}_2$, and with the quadratic Hamiltonians $\mathscr H_i$.
\end{proposition}

\begin{mproof}
We can use \cref{commzero} to compute
\begin{equation}\begin{split}
    [Q_m^{\gamma},Q_n^{\eta}] =& \int_{\gamma}\int_{\eta} \mathscr{P}(z)^{-m/2}\mathscr{P}(w)^{-n/2} [\varsigma_m(z)_{(0)},\varsigma_n(w)_{(0)}] \diff z \diff w\\
    = &\int_{\gamma}\int_{\eta} \mathscr{P}(z)^{-m/2}\mathscr{P}(w)^{-n/2} (\varsigma_m(z)_{(0)}\varsigma_n(w))_{(0)} \diff z \diff w,
\end{split}
\end{equation}
but we know from \cref{th2} that the 0th product between those states can be expressed as a sum of twisted derivatives and translations. It is now straightforward to check that the result of the commutator is zero: on one side because $(TX)_{(0)}=0$ for every state $X\in\mathbb{V}_0^{\bm{k}}$, on the other because of the property \eqref{Pfunction}.

To prove the second statement, we recall that $I^{a}_n = (I^{a}_{-1}\vakk)_{(n)}$ and the general property that for any two states $A,B$ of a vertex algebra we have
\begin{equation}
    [A_{(m)},B_{(n)}] = \sum_{k\geq 0} \binom{m}{k}(A_{(k)}B)_{(m+n-k)}, \qquad m,n\in\ZZ
\end{equation}
where
\begin{equation}
    \binom{m}{k} = \frac{m(m-1)\dots (m-k+1)}{k!}, \quad k\in\ZZ_{>0};\qquad \binom{m}{0}=1.
\end{equation}
Note that this represents a generalization of \cref{commzero}, see e.g. \cite{VY2017,frenkelvertex}.
At this point we can consider the generators $\{I_n^a\}_{a=1}^3$ with $n\in\ZZ$ and we get
\begin{equation}\begin{split}
[I^{a}_n,Q_m^{\gamma}] =&\int_{\gamma} \mathscr{P}(z)^{-m/2} [I^{a}_n,\varsigma_m(z)_{(0)}] \diff z\\
    =&\int_{\gamma}\sum_{k\geq0} \binom{n}{k}\mathscr{P}(z)^{-m/2} (I^{a}_k\varsigma_m(z))_{(n-k)} \diff z.
\end{split}
\end{equation}
It is now straightforward to check that the result is zero, by using property \eqref{singularity}, described in the relevant cases in \cref{singS1,singS3}, and the property \eqref{Pfunction}.
\end{mproof}

\subsubsection{Fourier modes}\label{fouriermodes}
Even though the operators $Q^{\gamma}_m$ we have just defined have all the right characteristics to be well-defined Hamiltonians as pointed out in \cref{prop}, there is one last subtlety about these objects, related to the fact that we want their action on highest weight modules to be diagonalisable.
In fact, considering $X^{(i)}_n\in \widetilde{U}(\widehat{\mathfrak{sl}}_2^{\oplus N})$ such that $\deg{X^{(i)}_n}=n$ and setting $\deg(\vakk)=0$, induces a $\ZZ_{\leq 0}$-gradation on the product of vacuum Verma modules, called the \textit{homogeneous gradation}. Therefore if $X\in\mathbb{V}_0^{\bm{k}}$ with degree $\deg(X)=k$, the degree of its modes is $\deg(X_{(m)})=1+k+m$.
The objects we have constructed $\varsigma_n(z)$, by definition, have $\deg(\varsigma_n(z))=-n-1$, therefore $\deg(\varsigma_n(z)_{(0)})=-n$.

This shows that in the homogeneous gradation these operators have degree $\neq 0$: this means that if we consider a module over $U(\widehat{\mathfrak{sl}}_2^{\oplus N})$ which has a trivial subspace of grade $n$ for large $n$, then the operator $\int_{\gamma_n} \mathscr{P}(z)^{-n/2}\varsigma_n(z)_{(0)}dz$ has a non-zero eigenvalue.

A way to overcome this issue is to consider the notion of Fourier mode $X_{[n]}\in \widetilde{U}(\widehat{\mathfrak{sl}}_2^{\oplus N})$ of the state $X\in\mathbb{V}_0^{\bm{k}}$: they have the property that we are looking for, namely $\deg(X_{[n]})=n$. Additionally, they satisfy a similar relation to \eqref{commzero},
\begin{equation}\label{fourierrel}
    [X_{[0]},Y_{[0]}]=(X_{(0)}Y)_{[0]},
\end{equation}
with $(T X)_{[0]}=0$. One has $(x^{(i)}_{-1}\vakk)_{[n]}=x_n$ for $x\in\sl_2$ and it is possible to show that the following recursive formula holds:
\begin{equation}\label{QI}
    (A_{(-n)}B)_{[m]}=((A\otimes f(t)) B)_{[m]}+\sum_{k>0}c_kA_{[-k]}B_{[k+m]}+\sum_{k\leq 0}c_kB_{[k+m]}A_{[-k]},
\end{equation}
where $f(t)$ is the Taylor series in $t:=u-v$ given by
\begin{equation}
    f=\frac{1}{(n-1)!}(-\partial_u)^{n-1}\Big(\frac{1}{u-v}-\iota_{u-v}\frac{e^{v}}{e^{u}-e^{v}}\Big).
\end{equation}
and where the coefficients $c_k$ are defined by the requirement that $\sum_{k>0}c_k(\frac{z}{w})^k$ and $-\sum_{k\leq0}c_k(\frac{z}{w})^k$ are the expansions, for $|z|<|w|$ and $|w|<|z|$ respectively, of the function
\begin{equation}
    \frac{1}{(n-1)!}(-w\partial_w)^{n-1}\frac{z}{w-z}.
\end{equation}
 The first relevant examples are
\begin{subequations}
\begin{align}
    \label{QI1}
    &(A_{(-1)}B)_{[m]}=
    \begin{aligned}[t]
         \frac{1}{2}(A_{(0)}B)_{[m]}-\frac{1}{12}(A_{(1)}&B)_{[m]}+\dots\\ &+\sum_{k>0}A_{[-k]}B_{[k+m]}+\sum_{k\leq 0}B_{[k+m]}A_{[-k]}
    \end{aligned}\\
    \label{QI2}
    &(A_{(-2)}B)_{[m]}=
    \begin{aligned}[t]
        \frac{1}{12}(A_{(0)}B)_{[m]}-\frac{1}{240}&(A_{(2)}B)_{[m]}+\dots\\ +&\sum_{k>0}kA_{[-k]}B_{[k+m]}+\sum_{k\leq 0}(-k)B_{[k+m]}A_{[-k]}
    \end{aligned}
\end{align}
\end{subequations}
where $A,B\in\mathbb{V}_0^{\bm{k}}$. These formulae are the Fourier-analog of the normal ordered product formula \cref{normalord}, and they allow one to compute by recursion the Fourier modes of a general state $X\in\mathbb{V}_0^{\bm{k}}$.

Property \eqref{fourierrel} means that if the vertex algebra 0th product of $X$ and $Y$ vanishes their Fourier zero-modes generate a commutative subalgebra, in homogeneous degree zero, of $\widetilde{U}(\widehat{\mathfrak{sl}}_2^{\oplus N})$. We let
\begin{equation}
   \widehat Q_n^{\gamma} = \int_{\gamma} \mathscr{P}(z)^{-n/2} \varsigma_n(z)_{[0]} \diff z ,
\end{equation}
for $n=1,3$. By the same logic as for \cref{prop}, we have the following.
\begin{proposition}\label{gaudinham}
The operators $\widehat Q_n^{\gamma}\in\widetilde{U}_{\bm{k}}(\hf{sl}_2^{\oplus N})$ have homogeneous degree 0 and they commute amongst themselves, with the generators of $\hf{sl}_2$, and with the quadratic Hamiltonians $\mathscr H_i$.
\end{proposition}

\section{Higher local Hamiltonians to sub-leading order}\label{subleading}
In the previous section, we have shown that it is possible to define quartic local Hamiltonians which commute among themselves and with the quadratic ones, together with the generators of $\sl_2$.
Following the same steps, one could in principle try to construct the Hamiltonians for every exponent of $\widehat{\mathfrak{sl}}_2$. However, the direct calculation (already lengthy in the case of $\varsigma_3(z)_{(0)}\varsigma_3(w)$, as we noted above) becomes computationally very demanding. What we shall do instead is work to next-to-leading order in a certain semiclassical limit, which will at least give a strong consistency check on the existence of the expected Hamiltonian densities.

Thus, let us introduce a formal parameter $\hbar$ and work over $\CC[[\hbar]]$. So in particular, all vector spaces above are now to be regarded as modules over $\CC[[\hbar]]$. We consider the following rescaled generators:
\begin{align}\label{classlim}
    I^a \longrightarrow \I^a &:= \hbar I^a, \nonumber\\
    \mathsf{k}   \longrightarrow \kk &:= \mathsf{k}.
\end{align}
With this re-scaling the commutation relations become
\begin{equation}
    [\I^{a \der p}_m(z),\I^{b \der q}_n(z)] = - \frac{p!q!}{(p+q+1)!} ( \hbar f^{ab}_c \I^{c\der{p+q+1}}_{m+n}(z) - \hbar^2 n \delta^{ab}\delta_{m+n,0}\kk^{\der{p+q+1}}(z)).
\end{equation}
At this point we can identify the various quantum corrections by their $\hbar$ dependence and work grade by grade.
We shall work at next-to-leading order, i.e. the next order beyond the usual semi-classical calculation of Poisson brackets. Thus, we consider the densities of Hamiltonians up to and including the leading quantum corrections at order $\hbar$, and we compute commutators up to and including terms of order $\hbar^2$.

\begin{remark*}
It is worth remarking that the classical limit, \cref{classlim}, that we take is not quite the standard one which recovers the usual classical Gaudin model (cf. \cite{LKT22} for a very complete discussion of that limit), because for us the central charges remain $O(1)$ in the limit.
From our present perspective this is simply for computational convenience -- this limit produces the simplest possible non-trivial check, and had we rescaled the central charges there would be more potential quantum correction terms already at next-to-leading order. But it might be interesting to consider this classical limit in its own right.
\end{remark*}
\medskip

Having introduced the formal parameter $\hbar$, there is a gradation on the enveloping algebras in which $\I^a$ and $\hbar$ have grade one and $\kk$ has grade zero.
Recall that $\mathcal V_n$ denotes the space of homogeneous meromorphic states of degree $n$, \cref{gradation}.
Let now $\widetilde{\mathcal{V}}_n\subset \mathcal V_n$ denote the subspace consisting of states that are also of grade $n$ in this new gradation (\emph{i.e.} which are sums of terms having exactly $n$ factors of $\I^a$ or $\hbar$).

\begin{proposition}
Modulo terms of order $\hbar^2$ there is, up to rescaling, exactly one state $\tilde{\varsigma}_{2n-1}\in\widetilde{\mathcal{V}}_{2n}^{\mathfrak{sl}_2}$, $n\in\ZZ_{\geq 1}$, such that, for all $x\in\sl_2$ and $m\in\ZZ_{\geq 0}$,
\begin{equation}
        \Delta x_m \tilde{\varsigma}_{2n-1} = 0 \qquad \textup{mod}\:\hbar^2 D^{(2n-1)}_z\mathcal{V}_{2n-m,2n-1}, \quad \textup{mod}\:\hbar^3\mathcal{V}_{2n-m,2n}.
\end{equation}
Explicitly, modulo terms in $\hbar^2\mathcal{V}_{2n}$,
\begin{equation}\label{uniquegen}\begin{split}
    \tilde{\varsigma}&_{2n-1}(z)= t_{i_1,\dots,i_{2n}} \I^{i_1}_{-1}(z)\I^{i_2}_{-1}(z)\cdots \I^{i_{2n}}_{-1}(z)\vac\\
    &\qquad+ \hbar \frac{n(2n+1)(2n-2)}{(2n-1)} t_{i_1,\dots,i_{2n-4}}f^{abc} \I^{a}_{-2}(z)\I^{b\prime}_{-1}(z)\I^{(c}_{-1}(z)\I^{i_1}_{-1}(z)\cdots \I^{i_{2n-4})}_{-1}(z)\vac.
\end{split}
\end{equation}
\end{proposition}

\vspace{.1cm}

\begin{mproof}
Given the basis $\{I^r\}_{r=1}^3$ for $\sl_2$, we need to show that there exist a function $\mathsf{G}_{m}(z)$ such that
\begin{equation}\label{singhbar}
    \Delta I^r_{m} \tilde{\varsigma}_{2n-1}(z) = D^{(2n-1)}_z \mathsf{G}^r_{m}(z) + \mathcal O(\hbar^3).
\end{equation}
For $m=0$, this is always true thanks to the invariance of the tensor \eqref{invariance2}.
For the same reasons explained in the previous section, the only relevant check that one needs to make is the one for $m=1$.
From direct calculation, we get
\begin{equation}
   \mathsf{G}^r_{1}(z) = -\frac{4n}{(2n-1)} \hbar^2 t_{i_1,\dots,i_{2n-2}} \I^{(r}_{-1}(z)\I^{i_1}_{-1}(z)\cdots \I^{i_{2n-2})}_{-1}(z)\vac.
\end{equation}
Note that this result is in accordance with the exact ones obtained for the quadratic ($n=1$) and the leading order of quartic ($n=2$) states, cf. \cref{G11,G13}.
\end{mproof}

(Observe that this is consistent with \cref{singS3} because the vectors in \cref{vectrans,vecTD} all come with factors of $\hbar^2$ in the limit.)

We can now state the following theorem
\begin{theorem}
Let $\tilde\varsigma_{n}(z)$ be as in \cref{uniquegen} above, for all odd $m,n\in\ZZ_{\geq1}$. We have
\begin{equation}\label{commhbar}
\begin{split}
 \tilde\varsigma_m(z)_{(0)}\tilde\varsigma_n(w) = (n D^{(m)}_z - m D^{(n)}_w)\mathsf{A}_{m,n}(z,w) + T\mathsf{B}_{m,n}(z,w)+\mathcal{O}(\hbar^3),
\end{split}
\end{equation}
 where $\mathsf{A}_{m,n}(z,w),\mathsf{B}_{m,n}(z,w)\in \widetilde{\mathcal{V}}_{(z,w)}$ are given by
\begin{equation}
    \begin{split}
       \mathsf{A}_{m,n}(z,w) = &\zeta_{m,n} \frac{\hbar^2}{z-w}t_{i_1,\dots,i_{m-1}}t_{j_1,\dots,j_{n-1}}\\
       &\times T\Big(\I_{-1}^{(a }(z)\I_{-1}^{i_1}(z)\dots \I_{-1}^{i_{m-1})}(z)\Big) \I_{-1}^{(a }(w)\I_{-1}^{j_1}(w)\dots \I_{-1}^{j_{n-1})}(w)\vac + \mathcal{O}(\hbar^3),
    \end{split}
\end{equation}
\begin{equation}
    \begin{split}
       \mathsf{B}_{m,n}(z,w) = &\zeta_{m,n} \frac{\hbar^2}{(z-w)^2}t_{i_1,\dots,i_{m-1}}t_{j_1,\dots,j_{n-1}}\\
       &\times \I_{-1}^{(a }(z)\I_{-1}^{i_1}(z)\dots \I_{-1}^{i_{m-1})}(z) \I_{-1}^{(a }(w)\I_{-1}^{j_1}(w)\dots \I_{-1}^{j_{n-1})}(w)\vac+\mathcal{O}(\hbar^3),
    \end{split}
\end{equation}
where
\begin{equation}
        \zeta_{m,n} = \frac{2(m+1)(n+1)}{mn}.
\end{equation}
Moreover, $\mathsf{A}_{m,n}(z,w)$ is a regular function for $z=w$ modulo translates and modulo terms proportional to $\hbar^3$.
\end{theorem}

\begin{mproof}
The 0th mode of $\varsigma_{m}(z)\in\widetilde{\mathcal{V}}_{m+1}$ can be inferred from a purely combinatorial reasoning.
Let us start with the top term $\tilde{\varsigma}_m^{\textup{TT}}(z)$ of \eqref{uniquegen}.
We know that computing the 0th mode, the number of generators in any term we get does not change, but the result will be a state of total depth $m$ and therefore we know there must be at least one generator with a positive mode. We can also use the fact that we are working at leading order in $\hbar$, therefore we could get at least one $\I_{0}$, one $\I_{1}$ or a term with two $\I_{0}$, every other term will be $\mathcal{O}(\hbar^3)$.
The only thing to fix is the combinatorial factor describing the number of possible ways to write such terms. The result is
\begin{equation*}
\begin{split}
    \tilde{\varsigma}_m^{\textup{TT}}(z)_{(0)}=t_{i_1,\dots,i_{m+1}}\Big[&\frac{(m+1)!}{(m-1)!}\I^{i_1}_{-2}(z)\I^{i_2}_{-1}(z)\dots \I^{i_{m}}_{-1}(z)\I^{i_{m+1}}_{1}(z)\\
        & +\frac{(m+1)!}{(m)!}\I^{i_1}_{-1}(z)\dots \I^{i_{m}}_{-1}(z)\I^{i_{m+1}}_{0}(z)\\
        &+ \frac{(m+1)!(m-1)}{2(m-1)!}\I^{i_1}_{-2}(z)\I^{i_2}_{-1}(z)\dots \I^{i_{m-1}}_{-1}(z)\I^{i_{m}}_{0}(z)\I^{i_{m+1}}_{0}(z)\Big] + \mathcal{O}(\hbar^3).
\end{split}
\end{equation*}
With similar arguments we can compute the 0th mode of the correction term $\tilde{\varsigma}_m^{\textup{C}}(z)$ of \eqref{uniquegen}, the result reads
\begin{equation*}
\begin{split}
    \tilde{\varsigma}_m^{\textup{C}}(z)_{(0)} =
    &-\hbar \xi (m-2)! t_{i_1,\dots,i_{m-3}} f^{abc}\I_{-2}^{a}(z) \I_{-1}^{(b }(z)\I_{-1}^{i_1}(z)\dots \I_{-1}^{i_{m-3})}(z)\I^c_0(z)\\
    &+\hbar \xi 2(m-2)(m-3)! t_{i_1,\dots,i_{m-3}} f^{abc}\I_{-2}^{a}(z)\I_{-1}^{b\prime}(z)\I_{-1}^{(i_1 }(z)\dots \I_{-1}^{i_{m-3})}(z)\I_{0}^{c}(z)\\
    &+\hbar \xi (m-2)(m-3)(m-3)! t_{i_1,\dots,i_{m-5}}f^{abc}\I_{-2}^{a}(z)\I_{-1}^{b\prime}(z)\\&\hspace{4.5cm}\times \I_{-1}^{(c}(z)\I_{-1}^{d }(z)\I^{i_1}_{-1}(z)\dots \I_{-1}^{i_{m-5})}\I_{0}^{d}(z)\\
    &-\hbar \xi  (m-2)(m-3)(m-3)!  t_{i_1,\dots,i_{m-5}} f^{abc}\I_{-2}^{d}(z)\I_{-1}^{a\prime}(z)\\
    &\hspace{4.5cm}\times \I_{-1}^{(b}(z)\I_{-1}^{d}(z)\I^{i_1}_{-1}(z)\dots \I_{-1}^{i_{m-5})}(z)\I_{0}^{c}(z) + \mathcal{O}(\hbar^3)\\
\end{split}
\end{equation*}
where $\xi = \frac{(m+2)(m+1)(m-1)}{2m}$. At this point, acting with what we have obtained on $\varsigma_{n}(w)$ and using repeatedly the commutation relations \eqref{commpar}, we obtain \cref{commhbar}.
\end{mproof}

\newpage
\appendix
\section{Full expression for \texorpdfstring{$\mathsf{A}_{1,3}(z,w)$}{TEXT} and \texorpdfstring{$\mathsf{B}_{1,3}(z,w)$}{TEXT} }\label{AppA}
The explicit expressions for $\mathsf{A}_{1,3}(z,w)$ and $\mathsf{B}_{1,3}(z,w)$ in $\mathcal V_{(z,w)}$ obtained by direct calculations are
\begin{equation*}
\begin{split}
\mathsf{A}_{1,3}(z,w) &=  \Big[\frac{8}{3}  \frac{1}{z-w}  I^{a}_{-2}(z) I^{(a}_{-1}(w)I^{b}_{-1}(w)I^{b)}_{-1}(w)\\
&+f^{abc}\Big( \frac{80}{9}  \frac{1}{z-w}  I^{a}_{-2}(z)I^{b}_{-2}(w)I^{c\prime}_{-1}(w)
-\frac{80}{9}  \frac{1}{(z-w)^2}  I^{a}_{-2}(z)I^{b}_{-2}(w)I^{c}_{-1}(w)\\
&\hspace{2cm}+\frac{160}{9}  \frac{1}{z-w}  I^{a}_{-3}(z)I^{b\prime}_{-1}(w)I^{c}_{-1}(w)\Big)\\
&-\frac{80}{9}  \frac{1}{z-w}  I^{a}_{-2}(z)I^{a}_{-3}(w)\mathsf{k}^{\prime}(w)
-\frac{80}{3}  \frac{1}{z-w}  I^{a}_{-4}(z)I^{a}_{-1}(w)\mathsf{k}^{\prime}(w)\\
&+\frac{320}{27}  \frac{1}{(z-w)^3}  I^{a}_{-2}(z)I^{a}_{-3}(w)
+\frac{160}{9}  \frac{1}{z-w}  I^{a}_{-4}(z)I^{a\prime\prime}_{-1}(w)\\
&+\frac{320}{9}  \frac{1}{(z-w)^3}  I^{a}_{-4}(z)I^{a}_{-1}(w)
+\frac{160}{27}  \frac{1}{(z-w)^2}  I^{a}_{-2}(z)I^{a\prime}_{-3}(w)\\
&-\frac{160}{9} \frac{1}{z-w}  I^{a}_{-3}(z)I^{a\prime\prime}_{-2}(w)
-\frac{160}{9}  \frac{1}{(z-w)^2}  I^{a}_{-3}(z)I^{a\prime}_{-2}(w)\\
&-\frac{160}{9} \frac{1}{(z-w)^2}  I^{a}_{-4}(z)I^{a\prime}_{-1}(w)\Big]\vac.
\end{split}
\end{equation*}

\begin{equation*}
    \begin{split}
\mathsf{B}_{1,3}(z,w) = &\Big[\frac{8}{3}  \frac{1}{(z-w)^2} I^{a}_{-1}(z)I^{(a}_{-1}(w)I^{b}_{-1}(w)I^{b)}_{-1}(w)\\
&+f^{abc}\Big(
-\frac{160}{9}  \frac{1}{(z-w)^3}  I^{a}_{-1}(z)I^{b}_{-2}(w)I^{c}_{-1}(w)
-\frac{80}{9}   \frac{1}{(z-w)^2}  I^{a}_{-1}(z)I^{b\prime}_{-1}(w)I^{c}_{-2}(w)\\
&\hspace{2cm}-\frac{160}{9}   \frac{1}{(z-w)^2}  I^{a}_{-2}(z)I^{b\prime}_{-1}(w)I^{c}_{-1}(w)\Big)\\
&-\frac{160}{9}  \frac{1}{(z-w)^2}  I^{a}_{-1}(z)I^{a}_{-3}(w) \mathsf{k}^{\prime}(w)
+\frac{1120}{27}   \frac{1}{(z-w)^2}  I^{a}_{-1}(z)I^{a\prime\prime}_{-3}(w)\\
&+\frac{640}{27}   \frac{1}{(z-w)^3}  I^{a}_{-1}(z)I^{a\prime}_{-3}(w)
-\frac{320}{9}   \frac{1}{(z-w)^2}  I^{a}_{-2}(z)I^{a\prime}_{-2}(w) \\
&+\frac{320}{9} \frac{1}{(z-w)^4}  I^{a}_{-1}(z)I^{a}_{-3}(w)
-\frac{640}{9} \frac{1}{(z-w)^3}  I^{a}_{-3}(z)I^{a\prime}_{-1}(w) \\
&+\frac{320}{3}\frac{1}{(z-w)^4}  I^{a}_{-3}(z)I^{a}_{-1}(w)\Big]\vac.
    \end{split}
\end{equation*}
\newpage
\printbibliography
\end{document}